\documentclass{amsart}
\usepackage{amsmath,amsthm}
\usepackage{amsfonts,amssymb}

 \usepackage{graphicx}
 \usepackage{ifpdf}

\newtheorem{thm}{Theorem}[section]

\newtheorem{lem}[thm]{Lemma}
\newtheorem{prop}[thm]{Proposition}

\newtheorem{defn}[thm]{Definition}

\theoremstyle{remark}

\makeatletter
\newcommand{\bddots}{%
  \mathinner{\mkern1mu\raise\p@\vbox{\kern7\p@\hbox{.}}\mkern2mu
    \raise4\p@\hbox{.}\mkern2mu\raise7\p@\hbox{.}\mkern1mu}}
\makeatother

\def\kb{{\mathbf \it k}}

\def\jb{{\mathbf j}}
\def\kb{{\mathbf k}}
\def\lb{{\mathbf l}}

\def\sb{{\mathbf s}}
\def\tb{{\mathbf t}}

 \def\CP{{\mathcal P}}
 
 \def\CI{{\mathcal I}}
 \def\CL{{\mathcal L}}
 \def\CT{{\mathcal T}}
 \def\CTC{{\mathcal T\!C}}
 \def\CTS{{\mathcal T\!S}}
 \def\CV{{\mathcal V}}
 \def\CH{{\mathcal H}}
 \def\CC{{\mathbb C}}
 \def\HH{{\mathbb H}}
 \def\PP{{\mathbb P}}
 \def\QQ{{\mathbb Q}}
 \def\RR{{\mathbb R}}
 \def\ZZ{{\mathbb Z}}
 \def\A{{\mathcal A}}
 \def\G{{\mathcal G}}
 \def\S{{\mathcal S}}
\newcommand{\e}{\mathrm{e}}
\newcommand{\ve}{\mathrm{v}}
\newcommand{\tr}{{\mathsf {tr}}}
\newcommand{\TC}{{\mathsf {TC}}}
\newcommand{\TS}{{\mathsf {TS}}}

 \def\diag{\operatorname{diag}}
 \def\sspan{\operatorname{span}}
\def \la {\langle}
\def \ra {\rangle}
\newcommand{\wt}{\widetilde}
\newcommand{\wh}{\widehat}

 \ifpdf
 \DeclareGraphicsExtensions{.pdf,.png,.jpg}
 \else
 \DeclareGraphicsExtensions{.eps}
 \fi
 \graphicspath{{images/}}

\begin{document}

\title[Discrete Fourier analysis on hexagon and triangle]
{Discrete Fourier analysis, Cubature and Interpolation on a Hexagon 
and a Triangle}
\author{Huiyuan Li} 
\address{Institute of Software\\
Chinese Academy of Sciences\\ Beijing 100080,China}
\email{hynli@mail.rdcps.ac.cn}
\author{Jiachang Sun}
\address{Institute of Software\\
Chinese Academy of Sciences\\ Beijing 100080,China}
\email{sun@mail.rdcps.ac.cn}
\author{Yuan Xu}
\address{  Department of Mathematics\\ University of Oregon\\
    Eugene, Oregon 97403-1222.}
\email{yuan@math.uoregon.edu}
 
\date{\today}
\thanks{The first and the second authors were supported by NSFC Grant 
10431050 and 60573023. The second author was supported by National Basic
Research Program Grant 2005CB321702. The third author was supported by 
NSF Grant DMS-0604056}

\begin{abstract}
Several problems of trigonometric approximation on a hexagon and a 
triangle are studied using the discrete Fourier transform and orthogonal
polynomials of two variables. A discrete Fourier analysis on the regular 
hexagon is developed in detail, from which the analysis on the triangle
is deduced. The results include cubature formulas and interpolation on these 
domains. In particular, a trigonometric Lagrange interpolation on a triangle 
is shown to satisfy an explicit compact formula, which is equivalent to 
the polynomial interpolation on a planer region bounded by Steiner's hypocycloid. 
The Lebesgue constant of the interpolation is shown to be in the order of
$(\log n)^2$. Furthermore, a Gauss cubature is established on the hypocycloid. 
\end{abstract}

\maketitle

\section{Introduction}
\setcounter{equation}{0}

Approximation by trigonometric polynomials, or equivalently complex 
exponentials, is at the root of approximation theory. Many problems 
central to approximation theory were first studied for trigonometric 
polynomials on the real line or on the unit circle, as illustrated in the 
classical treatise of Zygmund \cite{Z}. Much of the advance in the 
theory of trigonometric approximation is due to the periodicity of the function. 
In the case of one variable, one usually works with $2\pi$ periodic functions. 
A straightforward extension to several variables is the tensor product type,
where one works with functions that are $2\pi$ periodic in each of their 
variables. There are, however, other definitions of periodicity in several 
variables, most notably the periodicity defined by a lattice, which is a discrete
subgroup defined by $A\ZZ^d$, where $A$ is a nonsingular matrix, and
the periodic function satisifes $f(x+ A k) = f(x)$,  for all $k \in\ZZ^d$.
With such a periodicity, one works with exponentials or trigonometric functions of 
the form $\e^{i 2\pi \alpha \cdot x}$, where $\alpha$ and $x$ are in proper
sets of $\RR^d$, not necessarily the usual trigonometric polynomials. 

If $\Omega$ is a bounded open set that tiles $\RR^d$ with the lattice 
$L= A\ZZ^d$, meaning that $\Omega + A \ZZ^d = \RR^d$,  then a theorem of 
Fuglede \cite{F} states that the family of exponentials 
$\{\e^{ 2\pi i \alpha \cdot x} : \alpha \in L^\perp\}$, where $L^\perp = 
A^{-\tr}\ZZ^d$ is the dual lattice of $L$, forms an orthonormal 
basis in $L^2(\Omega)$. In particular, this shows that it is appropriate to use 
such exponentials for approximation. Naturally one can study Fourier series 
on $\Omega$, defined by these exponentials, which is more or less the 
usual multivariate Fourier series with a change of variables from $\ZZ^d$ to
$A\ZZ^d$.  Lattices appear prominently in various problems in 
combinatorics and analysis (see, for example, \cite{CS, E}), the Fourier 
series associated  with them are studied in the context of information
science and physics, such as sampling theory \cite{Hi,Ma} and multivariate 
signal processing \cite{DM}.  As far as we are aware, they have not been 
studied much in approximation theory until recently. Part of the reason may 
lie in the complexity of the set $\Omega$.  In the simplest non-tensor 
product case, it is a regular hexagon region on the plane, not one of those 
regions that people are used to in classical analysis. It turns out, 
however, that the results on the hexagon are instrumental for the study 
of trigonometric functions on the triangle \cite{Sun, LS}. The Fourier series
on the hexagon were studied from scratch in \cite{Sun}, and two families of 
functions analogues to those of sine and cosine but defined using three 
directions were identified and studied. In \cite{LS} the two generalized 
trigonometric functions were derived as solutions of the eigenfunctions 
of the Laplace equation on a triangle, and the Fourier series in terms of 
them were studied. 

The purpose of the present paper is to study discrete Fourier analysis
on the hexagon and on the triangle; in particular, we will derive results on 
cubature and Lagrange interpolation on these domains.  We shall follow 
two approaches. The first one relies on the discrete Fourier transforms 
associated with the lattices, which leads naturally to an interpolation 
formula, and the result applies to a general lattice. We then apply the 
result to the hexagonal tiling. The symmetry group of the hexagon lattice 
is the reflection group $\A_2$.  The generalized cosine and sine are 
invariant and anti-invariant projections, respectively, of the basic exponential 
function under $\A_2$, and they can be studied as functions defined 
on a fundamental domain, which is one of the equilateral triangles inside 
the hexagon. This allows us to pass from the hexagon to the triangle. 
The discrete Fourier transform can be regarded as a cubature formula.
In the case of a triangle, the trigonometric functions for which such a
cubature holds is a graded algebra due to a connection to polynomials 
of two variables, which shows that the cubature is a Gaussian type. 

The second approach uses orthogonal polynomials in two variables. 
It turns out that the generalized cosine and sine functions were studied 
earlier by Koornwinder in \cite{K}, who also started with eigenfunctions
of the Laplace operator and identified them as orthogonal polynomials on 
the region bounded by SteinerÕs hypocycloid, which he called Chebyshev 
polynomials of the first and the second kind, respectively. It is well known 
that orthogonal polynomials are closely related to cubature formulas and 
polynomial interpolation (see, for example, \cite{DX, My, St, X97, X00}). 
Roughly speaking, a cubature formula of high precision exists if the set of 
its nodes is an algebraic variety of a polynomial ideal generated by 
(quasi)orthogonal polynomials, and the size of the variety is equal to the 
codimension of the ideal. Furthermore, such a cubature formula arises 
from an interpolation polynomial. The rich structure of the generalized 
Chebyshev polynomials allow us to study their common zeros and establish 
the cubature and the polynomial interpolation. One particular result shows 
that Gaussian cubature formula exists for one integral on the region bounded 
by the hypocycloid. It should be mentioned that, unlike the case of one variable, 
Gaussian cubature formulas rarely exist in several variables.  

The most concrete result, in our view, is the trigonometric interpolation on 
the triangle. For a standard triangle $\Delta:=\{(x,y): x \ge 0, y\ge 0, x+y \le 1\}$, 
the set of interpolation points are $X_n:=\{(\tfrac{i}{n},  \tfrac{j}{n}): 0 \le i \le 
j \le n\}$ in $\Delta$. We found a compact formula for the fundamental 
interpolation functions in terms of the elementary trigonometric functions, 
which ensures that the interpolation can be computed efficiently. Moreover,  
the uniform operator norm (the Lebesgue constant) of the interpolation is shown 
to grow in the order of $(\log n)^2$ as $n$ goes to infinity. This is in sharp
contrast to the algebraic polynomial interpolation on the set $X_n$, which 
also exists and can be given by simple formulas, but has an undesirable 
convergence behavior (see, for example, \cite{R}). In fact, it is well known 
that equally spaced points are not good for algebraic polynomial interpolation. 
For example, a classical result shows that the polynomial interpolation of 
$f(x) = |x|$ on equally spaced points of $[-1,1]$ diverges at all points other
than $-1,0,1$. 

The paper is organized as follows. In the next section we state the basic 
results on lattices and discrete Fourier transform, and establish the connection
to the interpolation formula. Section 3 contains results on the regular hexagon, 
obtained from applying the general result of the previous section. Section 4 
studies the generalized sine and cosine functions, and the discrete Fourier
analysis on the triangle. In particular, it contains the results on the 
trigonometric interpolation. The generalized Chebyshev polynomials and 
their connection to the orthogonal polynomials of two variables are discussed
in Section 5, and used to establish the existence of Gaussian cubature, and 
an alternative way to derive interpolation polynomial.

\section{Continuous and discrete Fourier analysis with lattice}
\setcounter{equation}{0}

\subsection{Lattice and Fourier series} 
A lattice in $\RR^d$ is a discrete subgroup that contains $d$ linearly 
independent vectors,
\begin{align} \label{lattice} 
   L:= \left\{k_1 a_1 + k_2 a_2 +\cdots + k_d a_d:\ \ k_i \in \ZZ,\, i=1,2,\ldots,d\right\},
 \end{align}
where $a_1, \ldots, a_d$ are linearly independent column vectors in $\RR^d$. 
Let $A$ be the $d\times d$ matrix whose columns are $a_1,\ldots,a_d$. Then
$A$ is called a generator matrix of the lattice $L$. We can write $L$ as $L_A$ 
and a short notation for $L_A$ is $A \ZZ^d$;  that is, 
$$ 
L_A = A\ZZ^d = \left\{A k : \\  k \in \ZZ^d \right\}. 
$$
Any $d$-dimensional lattice has a dual lattice $L^\perp$ given by 
$$
    L^\perp = \{x \in \RR^d :  x \cdot y \in \ZZ \,\, \hbox{for all $y \in L$} \},
$$
where $x \cdot y$ denote the usual Euclidean inner product of $x$ and $y$. 
The generator matrix of the dual lattice $L_A^\perp$ is $A^{-\tr}: = (A^\tr)^{-1}$,
where $A^\tr$ denote the matrix transpose of $A$ (\cite{CS}). Throughout this paper
whenever we need to treat $x \in \RR^d$ as a vector, we treat it as a column 
vector,  so that $x \cdot y = x^{\tr} y$. 

A bounded set $\Omega$ of $\RR^d$ is said to tile $\RR^d$ with the lattice 
$L$ if
$$
  \sum_{\alpha \in L} \chi_\Omega (x + \alpha) = 1, \qquad 
      \hbox{for almost all $x \in \RR^d$},
$$
where $\chi_\Omega$ denotes the characteristic function of $\Omega$. We 
write this as $\Omega + L = \RR^d$. Tiling and Fourier analysis are closely
related; see, for example, \cite{F,Ko}. We will need the following result due 
to Fuglede \cite{F}.  Let $\langle \cdot, \cdot\rangle_\Omega $ denote the
inner product in the space $L^2(\Omega)$ defined by 
\begin{equation} \label{ipOmega}
  \langle f, g \rangle_\Omega = \frac{1}{\mathrm{mes}(\Omega)} 
      \int_\Omega f(x) \overline{g(x)} dx. 
\end{equation}

\begin{thm} \cite{F}
Let $\Omega$ in $\RR^d$ be a bounded open domain and $L$ be a 
lattice of $\RR^d$. Then $\Omega + L = \RR^d$ if and only if 
$\{\e^{2 \pi i \alpha \cdot x}: \alpha \in L^\perp\}$ is an orthonormal 
basis with respect to the inner product \eqref{ipOmega}. 
\end{thm} 

The orthonormal property is defined with respect to normalized Lebesgue
measure on $\Omega$. If $L = L_A$, then the measure of $\Omega$ is 
equal to $|\det A|$. Furthermore, since $L^\perp_A = A^{-\tr}\ZZ^d$, we can
write $\alpha = A^{-\tr} k$ for $\alpha \in L_A^\perp$ and $k \in \ZZ^d$,
so that $\alpha \cdot x = (A^{-\tr} k) \cdot x = k^{\tr} A^{-1} x$. Hence the 
orthogonality in the theorem is 
\begin{equation} \label{expOrth}
  \frac{1}{|\det (A)|} \int_{\Omega} \e^{2\pi i k^{\tr} A^{-1} \xi}  d\xi =  \delta_{k,0}, 
      \qquad k \in \ZZ^d. 
\end{equation}
The set $\Omega$ is called a spectral (fundamental) set for the lattice $L$. 
If $L =L_A$ we also write $\Omega = \Omega_A$.

A function $f \in L^1(\Omega_A)$ can be expanded into a Fourier series 
$$
   f(x) \sim \sum_{k\in \ZZ^d} c_k \e^{2\pi i k^\tr A^{-1} x}, \qquad 
         c_k =   \frac{1}{|\det (A)|} \int_{\Omega} f(x) \e^{-2\pi i k^{\tr} A^{-1} x}  dx,
$$
which is just a change of variables away from the usual multiple Fourier series. 
The Fourier transform $\wh f$ of a function $f \in L^1(\RR^d)$ and its inverse 
are defined by 
\begin{align*}
\wh f(\xi) = \int_{\RR^d} f(x) \e^{-2\pi i\, \xi \cdot x} d x \quad\hbox{and}\quad
        f(x)    =  \int_{\RR^d} \wh f(\xi)  \e^{2\pi i\, \xi \cdot x} d \xi. 
\end{align*}
One important tool for us is the Poisson summation formula given by 
\begin{equation} \label{poisson}
    \sum_{k \in \ZZ^d} f (x+ A k) = \frac{1}{|\det (A)|} 
            \sum_{k \in \ZZ^d} \wh f(A^{-\tr} k) \e^{2\pi i k^\tr A^{-1} x}, 
\end{equation}
which holds for $f\in L^1(\Omega_A)$ under the usual assumption on the
convergence of the series in both sides. This formula has been used by 
various authors (see, for example, \cite{AB,E,Hi,K,SO}). An immediate
consequence of the Poisson summation formula is the following sampling 
theorem (see, for example, \cite{Hi,Ma}).

\begin{prop} \label{prop:sampling}
Let $\Omega$ be a spectral set of the lattice $A\ZZ^d$. Assume that
$\wh f$ is supported on $\Omega$ and $\wh f  \in L^2(\Omega)$. Then 
\begin{equation} \label{sampling}
      f(x) = \sum_{k \in \ZZ^d} f(A^{-\tr} k) \Psi_\Omega(x- A^{- \tr}k)
\end{equation}
in $L^2(\Omega)$, where 
\begin{equation} \label{Psi}
 \Psi_\Omega(x) = \frac{1}{|\det (A)|} \int_{\Omega} \e^{2 \pi i x\cdot \xi} d \xi. 
\end{equation}
\end{prop} 

In fact, since $\wh f$ is supported on $\Omega_A$, $\wh f (\xi) = 
\sum_{k \in \ZZ^d} \wh f(\xi+ Ak)$. Integrating over $\Omega$,  the sampling 
equation follows from \eqref{poisson} and Fourier inversion.  We note that
$$
\Psi_\Omega(A^{-\tr}j) = \delta_{0,j}, \qquad \hbox{ for all $j \in \ZZ^d$},
$$ 
by \eqref{expOrth}, so that $\Psi_\Omega$ can be considered as a 
cardinal interpolation function. 


\subsection{Discrete Fourier analysis and interpolation}

A function $f$ defined on $\RR^d$ is called a periodic function with respect
to the lattice $A \ZZ^d$ if 
$$
      f (x + A k) = f(x) \qquad \hbox{for all $k \in \ZZ^d$}.
$$
For a given lattice $L$, its spectral set is evidently not unique. The orthogonality
\eqref{expOrth} is independent of the choice of $\Omega$. This can be seen
directly from the fact that the function $x \mapsto e^{2\pi i k^\tr A^{-1} x}$ is 
periodic with respect to the lattice $A\ZZ^d$. 

In order to carry out a discrete Fourier analysis with respect to a lattice, we 
shall fix $\Omega$ such that $\Omega$ contains $0$ in its interior and we 
further require in this subsection that the tiling $\Omega + A\ZZ^d = \RR^d$ 
holds pointwise and without overlapping. In other words, we require
\begin{equation} \label{Omega}
\sum_{k \in \ZZ^d} \chi_\Omega(x+ Ak) = 1, \,\, \forall x\in \RR^d,
\quad \hbox{and}\quad
(\Omega + Ak) \cap (\Omega+ Aj) = \emptyset, \,\, \hbox{if $k \ne j$}.
\end{equation}
For example, we can take $\Omega = [-\tfrac{1}{2}, \tfrac{1}{2})^d$ for the 
standard lattice $\ZZ^d$.

\begin{defn} \label{def:N}
Let $A$ and $B$ be two nonsingular matrices such that $\Omega_A$ and 
$\Omega_B$ satisfy \eqref{Omega}.  Assume that all entries of $N := B^\tr A$ 
are integers. Denote 
$$
\Lambda_N: = \{ k \in \ZZ^d: B^{-\tr}k \in \Omega_A\}, \quad \hbox{and}
\quad
\Lambda_{N^\tr}: = \{ k \in \ZZ^d: A^{-\tr}k \in \Omega_B\}. 
$$
\end{defn}

Two points $x, y \in \RR^d$ are said to be congruent with respect to 
the lattice $A\ZZ^d$, if $x - y \in A\ZZ^d$, and we write $x \equiv y 
\mod A$.  

\begin{thm} \label{prop:d-ortho}
Let $A, B$ and $N$ be as in Definition \ref{def:N}. Then 
\begin{equation}\label{d-ortho}
  \frac{1}{|\det(N)|}\sum_{j \in \Lambda_N} \e^{2 \pi i k^\tr N^{-1} j} 
      = \begin{cases} 1,  &  \hbox{if  $k \equiv 0 \mod N^\tr$}, \\
              0, & \hbox{otherwise}. \end{cases}
\end{equation}
\end{thm}

\begin{proof}
Let $f \in L^1(\RR^d)$ be such that Poisson summation formula holds.  
Denote
\begin{align*}
       \Phi^A f (x): = \sum_{k \in \ZZ^d} f (x+ A k). 
\end{align*}
Then the Poisson summation formula \eqref{poisson} shows that 
\begin{align*}
  \Phi^A f(B^{-\tr} j ) = \frac{1}{|\det(A)|} \sum_{k \in \ZZ^d}
         \wh f(A^{-\tr} k)  e^{2 \pi i k^\tr N^{-1} j}, \qquad j \in \ZZ^d.
\end{align*}
Since  $\RR^d = \Omega_B + B\ZZ^d$, we can write $A^{-\tr} k = x_B + B m$ 
for every $k \in \ZZ^d$ with $x_B \in \Omega_B$ and $m \in \ZZ^d$. The fact 
that $N$ has integer entries shows that $A^\tr x_B = l  \in \ZZ^d$. That is, 
every $k \in \ZZ^d$ can be written as $k = l +  N^\tr m$ with 
$l \in \Lambda_{N^\tr}$ and $m \in \ZZ^d$. It is easy to see that
$\e^{2 \pi i k^\tr N^{-1} j }  = \e^{2 \pi i l^\tr N^{-1} j } $ under this decomposition. Consequently,
\begin{align} \label{PhiA}
  \Phi^A f(B^{-\tr} j) & = \frac{1}{|\det(A)|} \sum_{l \in \Lambda_{N^\tr}}
   \sum_{m \in \ZZ^d} \wh f (A^{-\tr} l + B m) \e^{2\pi i l^\tr N^{-1} j} \\
    & = \frac{1}{|\det(A)|} \sum_{l \in \Lambda_{N^\tr}} \Phi^B \wh f(A^{-\tr} l)
          \e^{2 \pi i l^\tr N^{-1} j}.  \notag
\end{align}  
Similarly, we have 
\begin{align} \label{PhiB}
  \Phi^B \wh f(A^{-\tr} l ) = \frac{1}{|\det(B)|} 
        \sum_{j \in \Lambda_{N}} \Phi^A f(B^{-\tr} j)
                      \e^{- 2 \pi i l^\tr N^{-1} j}.  
\end{align}  
Substituting \eqref{PhiA} into \eqref{PhiB} leads to the identity
$$
  \Phi^B \wh f(A^{-\tr} l ) = \frac{1}{|\det(AB)|} \sum_{k\in \Lambda_N^{\tr}}
          \Phi^B \wh f (A^{-\tr} k) \sum_{j \in \Lambda_{N}} 
                      \e^{2 \pi i (k-l)^\tr N^{-1} j} \, , 
$$ 
that holds
for $f$ in, say, the Schwartz class of functions, from which \eqref{d-ortho}
follows immediately.
\end{proof}

\begin{thm} \label{prop:d-inner}
Let $A, B$ and $N$ be as in Definition \ref{def:N}.  Define
$$
    \langle f, g \rangle_N =  \frac{1}{|\det (N)|} 
        \sum_{j \in \Lambda_N } f(B^{-\tr} j ) \overline{g(B^{-\tr} j )}
$$
for $f, g$ in $C(\Omega_A)$, the space of continuous functions on $\Omega_A$. 
Then
\begin{equation}\label{c-d-inner}
    \langle f, g \rangle_{\Omega_A} =  \langle f, g \rangle_{N}
\end{equation}
for all $f, g$ in the finite dimensional subspace
$$
  \CH_N := \sspan \left\{\phi_k:   \phi_k(x) = \e^{2\pi i k^\tr A^{-1} x}, \, 
              k \in \Lambda_{N^\tr} \right \}.   
$$
\end{thm}

\begin{proof}
By the definition of $\langle \cdot, \cdot \rangle_N$, equation \eqref{d-ortho} 
implies that $\langle \phi_l, \phi_k\rangle_N = \delta_{l,k}$, $l,k \in N^\tr$, which 
shows, by \eqref{expOrth}, that \eqref{c-d-inner} holds for $f= \phi_l$ and
$g =\phi_k$. Since both inner products are sesquilinear, the stated result 
follows. 
\end{proof}

If $\Lambda$ is a subset of $\ZZ^d$, we denote by $\#\Lambda$ the number
of points in $\Lambda$. Evidently the dimension of $\CH_N$ is 
$\# \Lambda_{N^\tr}$. 

Let $ \CI_N f $ denote the Fourier expansion of $f \in C(\Omega_A)$ in 
$\CH_N$ with respect to the inner product $\langle \cdot, \cdot \rangle$. Then 
\begin{align*}
   \CI_N f(x): = & \sum_{j\in \Lambda_{N^\tr} }  \langle f, \phi_j \rangle_N  \phi_j (x) \\ 
           = & \frac{1}{|\det(N)|} \sum_{k \in \Lambda_N } f(B^{-\tr}k)
                \sum_{j\in \Lambda_{N^\tr}}  \phi_j (x) \overline{\phi_j(B^{-\tr}k)}.
\end{align*}
Hence, analogous to the sampling theorem in Proposition \ref{prop:sampling}, we 
have 
\begin{equation} \label{interpolation}
 \CI_N f (x) = \sum_{k \in \Lambda_N} f(B^{-\tr}k)\Psi_{\Omega_B}^A(x-B^{-\tr}k),
 \qquad f \in C(\Omega_A), 
\end{equation}
where 
\begin{equation} \label{ell}
  \Psi_{\Omega_B}^A (x) = \frac{1}{|\det (N)|} \sum_{j \in \Lambda_{N^\tr}} 
       \e^{2 \pi i  j^\tr A^{-1} x}. 
\end{equation}

\begin{thm} \label{prop:interpolation}
Let $A, B$ and $N$ be as in Definition \ref{def:N}. Then $\CI_N f$ is the 
unique interpolation operator on $N$ in $\CH_N$; that is, 
$$
    \CI_N f (B^{-\tr}j ) = f (B^{-\tr}j),\qquad  \forall j \in \Lambda_N. 
$$
In particular, $\# \Lambda_N = \# \Lambda_{N^\tr} = |\det(N)|$. 
Furthermore, the fundamental interpolation function $\Psi_{\Omega_A}^B$ 
satisfies
\begin{equation}\label{Psi_Psi}
  \Psi_{\Omega_B}^A(x) = \sum_{j \in \ZZ^d} \Psi_{\Omega_B}(x+A j). 
 \end{equation}
\end{thm}

\begin{proof}
Using \eqref{d-ortho} with $N$ in place of $N^\tr$ gives immediately that $\Psi_{\Omega_B}^A(B^{-\tr}k)  = \delta_{k,0}$ for $k \in \Lambda_N$, so that the 
interpolation holds. This also shows that $\{\Psi_{\Omega_B}^A(x- B^{-\tr} k)
: k \in \Lambda_N\}$  is linearly independent.
The equation \eqref{d-ortho} with $k =0$ shows immediately  that $\# \Lambda_N
= |\det (N)|$. Similarly, with $N$ replaced by $N^\tr$, we have 
 $\# \Lambda_N = |\det (N^\tr)|$. This shows that   $\# \Lambda_N = |\det (N)|
 = \# \Lambda_{N^\tr}$. 
Furthermore, let $M = ( \phi_{k}(B^{-\tr}j))_{j,k \in \Lambda_N}$ be the interpolation
matrix, then the equation \eqref{d-ortho} shows that $M$ is a unitary matrix,
which shows in particular that $M$ is invertible. Consequently the interpolation
by $\CH_n$ on points in $\Lambda_N$  is unique. 

Finally, since  $\wh \Psi_{\Omega_B} = |\det (B)|^{-1}\chi_{\Omega_B}$,  we have 
 \begin{align*}
   \Psi_{\Omega_B}^A(x) & = \frac{1}{|\det(N)|} 
        \sum_{j\in \ZZ^d}  \phi_j(x) \chi_{\Omega_B}(A^{-\tr}j) \\
        & = \frac{1}{|\det(A)|} 
           \sum_{j\in \ZZ^d}  \phi_j(x) \wh \Psi_{\Omega_B}(A^{-\tr}j),
\end{align*}
from which \eqref{Psi_Psi} follows from the Poisson summation formula
\eqref{poisson}. 
\end{proof} 

The result in this section applies to fairly general lattices. In the following 
section we apply it to the hexagon lattice in two dimension. We hope to 
report results for higher dimensional lattices in a future work.


\section{Discrete Fourier analysis on the hexagon}
\setcounter{equation}{0}

\subsection{Hexagon lattice and Fourier analysis}
The generator matrix and the spectral set of the hexagon lattice is given by 
$$
H=\begin{pmatrix} \sqrt{3} & 0\\ -1 & 2\end{pmatrix}, \qquad 
\Omega_H =\left\{(x_1,x_2):\  -1\leq x_2, \tfrac{\sqrt{3}}{2}x_1 \pm 
   \tfrac{1}{2} x_2 < 1 \right\}.
$$ 
The strict inequality in the definition of $\Omega_H$ comes from  our assumption 
in \eqref{Omega}. 

\begin{figure}[htb]
\centerline{\includegraphics[width=0.5\textwidth]{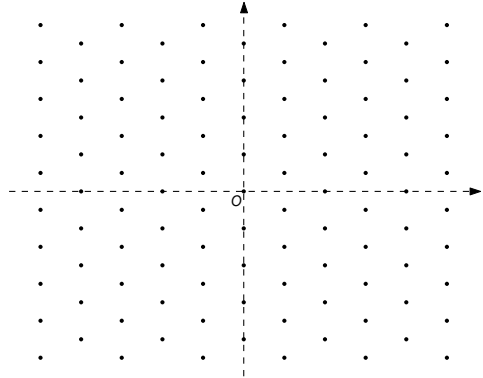}}
\caption{The lattice $L_{H}$.}
\end{figure}

As shown in \cite{Sun}, it is more convenient to use homogeneous coordinates 
$(t_1,t_2,t_3)$ that satisfies $t_1 + t_2 +t_3 =0$. We define  
\begin{align}\label{coordinates}
 t_1= -\frac{x_2}{2} +\frac{\sqrt{3}x_1}{2},\quad t_2 = x_2, \quad t_3 = 
      -\frac{x_2}{2} -\frac{\sqrt{3}x_1}{2}. 
\end{align}
Under these homogeneous coordinates, the hexagon $\Omega_H$ becomes 
\begin{align*}
\Omega=\left\{(t_1,t_2,t_3)\in \RR^3:\  -1\le  t_1,t_2,-t_3<1;\, 
    t_1+t_2+t_3=0 \right\},
\end{align*}
which is the intersection of the plane $t_1+t_2+t_3=0$ with the cube $[-1,1]^3$,
as shown in Figure 2. Later in the paper we shall depict the hexagon as a two
dimensional figure, but still refer to its points by their homogeneous coordinates,
as shown by the right hand side figure of Figure 2. 

\begin{figure}[h]
\centerline{\includegraphics[width=0.4\textwidth]{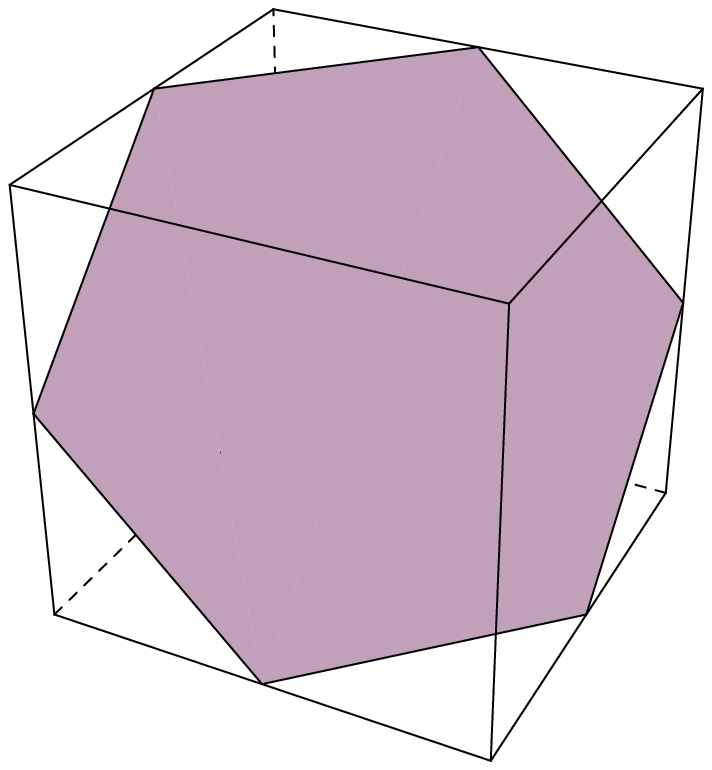} \quad
\includegraphics[width=0.5\textwidth]{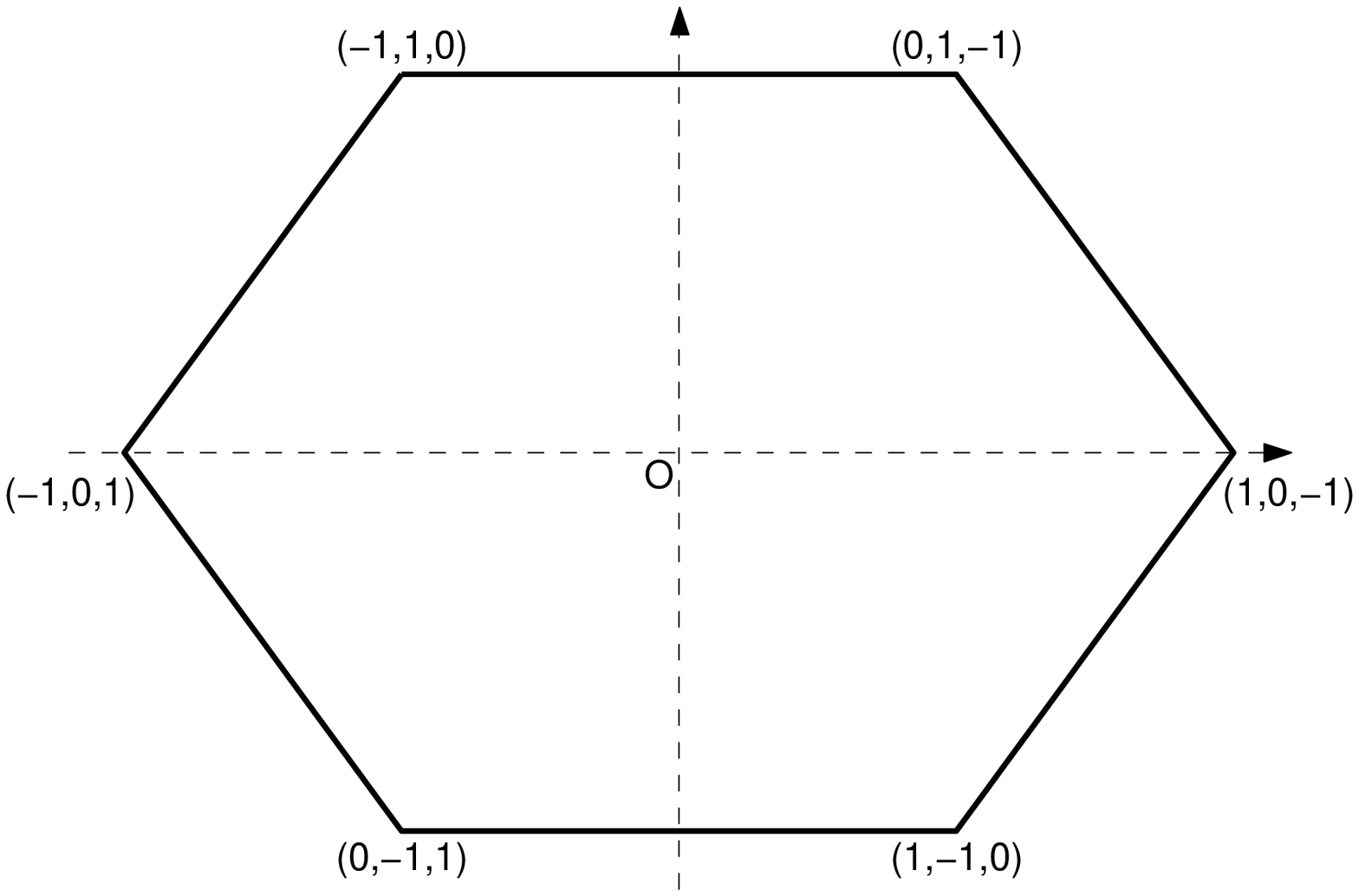}}
\caption{The hexagon in homogeneous coordinates.}
\end{figure}

For convenience, we adopt the convention of using bold letters, such 
as $\tb$, to denote points in the space 
$$
    \RR_H^3 : = \{\tb = (t_1,t_2,t_3)\in \RR^3: t_1+t_2 +t_3 =0\}. 
$$
In other words, bold letters such as  $\tb$ stand for homogeneous 
coordinates. Furthermore, we introduce the notation
$$
\HH := \ZZ^3 \cap \RR^3_H= 
     \left\{ \kb=(k_1,k_2,k_3)\in \ZZ^3:\ k_1 +k_2 +k_3 = 0\right\}.
$$
The inner product on the hexagon under homogeneous coordinates 
becomes 
\begin{align*}
\langle f, g\rangle_H = \frac{1}{|\Omega_H|} \int_{\Omega_H} f(x_1,x_2) 
   \overline{g(x_1,x_2)} d x_1 dx_2
= \frac{1}{|\Omega|} \int_{\Omega} f(\tb) \overline{g(\tb)} d \tb,
\end{align*}
where $|\Omega|$ denotes the area of $\Omega$. 

To apply the general results in the previous section, we choose $A = H$
and $B = \frac{n}{2}  H$, where $n$ is a positive integer. Then
$$
   N= B^\tr A = \frac{n}{2} H^T H = 
             \left[\begin{matrix} 2n &  -n \\ -n & 2n \end{matrix}\right]
$$
has integer entries. Note that $N$ is a symmetric matrix so that $\Lambda_N
= \Lambda_{N^\tr}$. Using the fact that 
$B^{-\tr} j = \tfrac{2}{n}H^{-\tr}j =  \tfrac{1}{n} \left( \tfrac{2}{\sqrt{3}}j_1 + 
\tfrac{1}{\sqrt{3}} j_2, j_2\right)$, it is not hard to see that $B^{-\tr} j \in \Omega_A$, 
or equivalently $j \in \Lambda_N$, becomes $\jb=(j_1,j_2,-j_1-j_2) \in \HH_n$,
where 
$$
    \HH_n : = \{\jb \in \HH: -n \le j_1,j_2, - j_3 < n \} \qquad \hbox{and}
      \qquad  \# \HH_n= 3n^2.
$$
In other words, $\Lambda_N$ in the previous section becomes $\HH_n$ 
in homogeneous coordinates. The fact that $\# \HH_n = 3 n^2$ can
be easily verified. Moreover, under the change of variables 
$x=(x_1,x_2) \mapsto  (t_1,t_2,t_3)$ in \eqref{coordinates}, it is easy to see that 
\begin{equation} \label{kAx}
      (k_1,k_2) H^{-1} x  =  \tfrac{1}{3} (k_1,k_2,-k_1-k_2) \tb 
             = \tfrac{1}{3} \kb \cdot \tb, 
\end{equation}
so that $\e^{2 \pi i k^{\tr} H^{-1} x} = \e^{\frac{2 \pi i}{3}\kb^\tr \tb}$. Therefore,
introducing  the notation 
$$
     \phi_\jb(\tb) : = \e^{\frac{2 \pi i}{3}\jb^\tr \tb}, \qquad \jb \in \HH,
$$
the orthogonality relation \eqref{expOrth} becomes 
\begin{equation}\label{H-expOrth}
    \langle\phi_\kb, \phi_\jb\rangle_H =\delta_{\kb,\jb}, \qquad \kb, \jb \in \HH. 
\end{equation}
The finite dimensional space $\CH_N$ in the previous section becomes
$$
   \CH_n: = \sspan \left \{ \phi_\jb: \jb \in \HH_n \right \}, \quad 
   \hbox{and} \quad \dim \CH_n = 3 n^2. 
$$

Under the homogeneous coordinates \eqref{coordinates},  
$x\equiv y \pmod{H}$ becomes  $\tb  \equiv \sb \mod 3$, where we define
$$
 \tb \equiv \sb \mod 3 \quad \Longleftrightarrow \quad  t_1-s_1 \equiv t_2-s_2 
      \equiv t_3-s_3 \mod 3.
$$
We call a function $f$ H-periodic if $f(\tb ) = f(\tb + \jb)$ whenever  
$\jb  \equiv 0 \mod 3$.  Since $\jb, \kb \in \HH$ implies that 
$2 \jb \cdot \kb = (j_1-j_2)(k_1-k_2) + 3 j_3k_3$, we see that 
$\phi_\jb$ is H-periodic. 
The following lemma is convenient in practice.

\begin{lem} \label{lem:H-periodic}
Let $\varepsilon_1 = (2,-1,-1)$, $\varepsilon_2 = (-1,2,-1)$, 
$\varepsilon_3 = (-1,-1,2)$.  Then a function $f(\tb)$ is H-periodic if and 
only if  
\begin{equation} \label{H-periodic}
     f(\tb + j \varepsilon_i) = f(\tb),  \quad j\in \ZZ, \quad i = 1, 2, 3. 
\end{equation}
\end{lem} 


\subsection{Discrete Fourier analysis on the regular hexagon} 

Using the set-up in the previous subsection, Theorem \ref{prop:d-inner}
becomes, in homogeneous coordinates,  the following: 

\begin{prop} \label{prop:H-inner}
For $n \ge 0$, define 
\begin{equation} \label{H-inner}
      \langle f, g\rangle_n := \frac{1}{3 n^2} \sum_{\jb \in \HH_n} f(\tfrac{\jb}{n})
           \overline{g(\tfrac{\jb}{n})}, \qquad f,g \in C(\Omega). 
\end{equation}
Then
$$
   \langle f, g\rangle_H =  \langle f, g\rangle_n,  \qquad f, g \in \CH_n.
$$
\end{prop}

The point set  $\HH_n$, hence the inner product $\langle \cdot, 
\cdot \rangle_n$, is not symmetric on $\Omega$ in the sense that it contains
only part of the points on the boundary, as shown by the left hand figure in 
Figure 3. We denote by $\HH^*_n$ the symmetric point set on $\Omega$, 
$$
    \HH_n^* : = \left\{\jb \in \HH: -n \le j_1,j_2,j_3 \le n \right\}, \quad
      \hbox{and} \quad \# \HH_n^* = 3n^2 + 3n +1.
$$ 
The set $\HH_n^*$ is shown by the right hand figure in Figure 3. 

\begin{figure}[h]
\centerline{\includegraphics[width=0.48\textwidth]{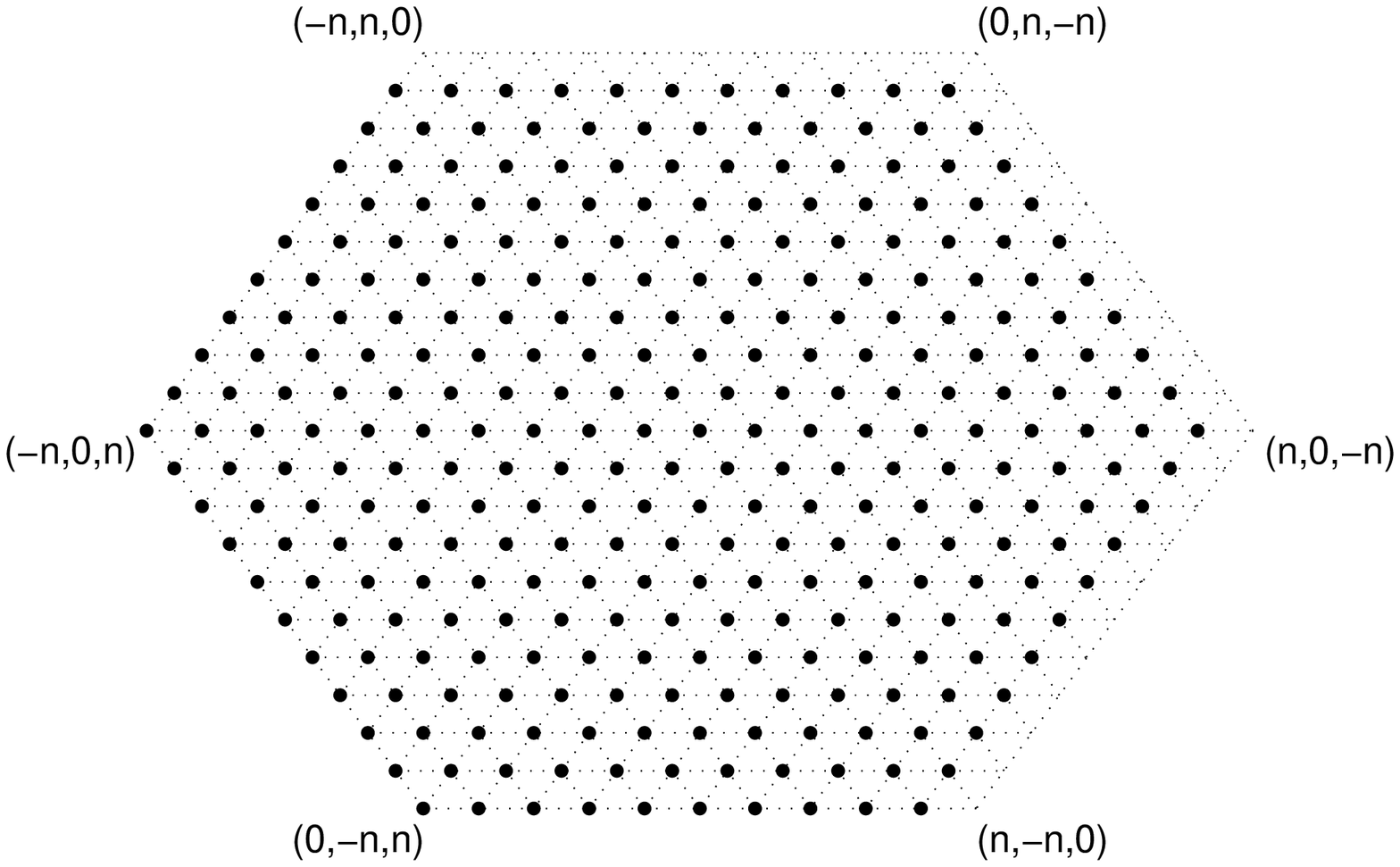}
 \quad \includegraphics[width=0.48\textwidth]{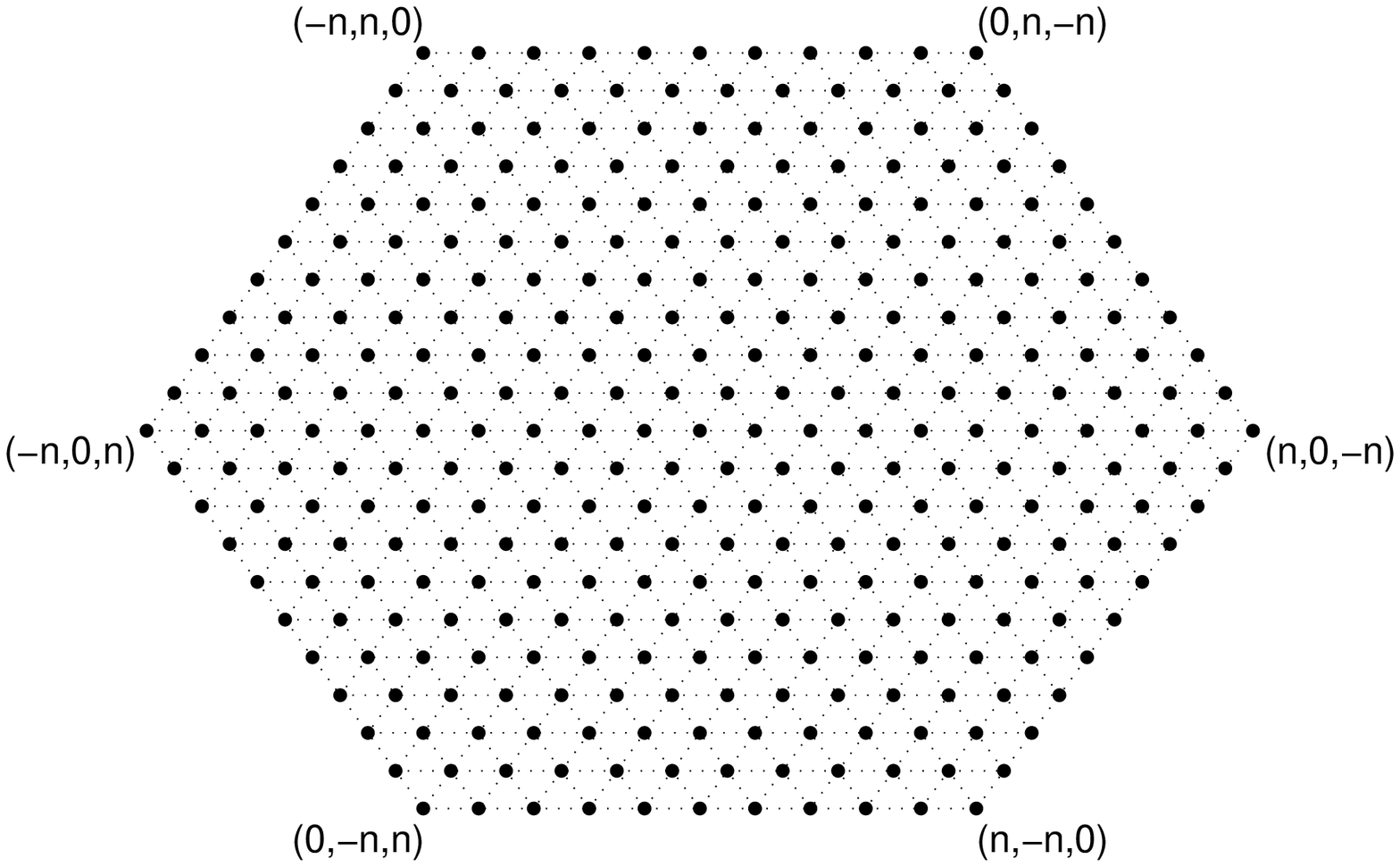}}
\caption{ The set $\HH_n$  and the set $\HH_n^*$.}
 \end{figure}

Using periodicity, we can show that the inner product $\langle \cdot, 
\cdot \rangle$ is equivalent to a symmetric discrete inner product based on 
$\HH_n^*$. Let 
$$
\HH_n^\circ : = \left\{\jb \in \HH: -n <  j_1,j_2,j_3 <  n \right\}
$$
denote the set of interior points in $\HH_n^*$ and $\HH_n$.  Define 
\begin{equation} \label{Symm-inner}
      \langle f, g\rangle_n^* := \frac{1}{3 n^2} \sum_{\jb \in \HH_n^*} 
           c_\jb^{(n)} f(\tfrac{\jb}{n})
           \overline{g(\tfrac{\jb}{n})}, \qquad f,g \in C(\Omega), 
\end{equation} 
where 
$$
c_{\jb}^{(n)} = \begin{cases} 1, & \jb \in \HH_n^\circ, \quad \text{(inner points)},\\
       \frac{1}{2},
        & \jb \in \HH_n^* \setminus \HH_n^\circ,\,  j_1j_2j_3\neq 0,\quad  
          \text{(edge points)},\\
       \frac{1}{3}, 
       & \jb \in \HH_n^* \setminus \HH_n^\circ, \, j_1j_2j_3 =  0,\quad  \text{(vertices)},
              \end{cases}
$$
in which we also include the positions of the points; the only one that needs
an explanation is the ``edge points", which are the points located on the edges 
but not on the vertices of the hexagon.  In analogy to Proposition \ref{prop:H-inner},
we have the following theorem.  

\begin{thm}\label{prop: Symm-inner}
For $n \ge 0$, 
$$
   \langle f, g\rangle_H = \langle f, g\rangle_n =  
    \langle f, g\rangle_n^*,  \qquad f, g \in \CH_n.
$$
\end{thm}

\begin{proof}
If $f$ is an H-periodic function, then it follows from \eqref{H-periodic} that  
\begin{align*}
     6 \sum_{\jb \in \HH^*_n\setminus \HH_{n}} & c_\jb^{(n)} f(\tfrac{\jb}{n})
       = 2 \left [ f(1,-1,0) + f(1,0,-1) + f(0,1,-1) + f(-1,1,0) \right ]  \\
      &\  + 3 \sum_{1\leq j \leq n-1} \left[
       f(1,- \tfrac{j}{n},\tfrac{j}{n}-1) + f(\tfrac{j}{n}-1,1,-\tfrac{j}{n}) +
               f(1-\tfrac{j}{n},\tfrac{j}{n},-1) \right]\\
       =&\  3 \sum_{1\leq j \leq n-1} \left[
       f(-1,1-\tfrac{j}{n},\tfrac{j}{n}) + f(\tfrac{j}{n},-1,1-\tfrac{j}{n}) + 
              f(-\tfrac{j}{n},\tfrac{j}{n}-1,1) \right]\\ 
       & \qquad \quad \ + 4 \left [ f(0,-1,1) + f(-1,0,1)\right] \\
       = & \   6 \sum_{\jb\in \HH_n\setminus \HH_{n-1}^*} 
                   \left (1 - c_\jb^{(n)} \right) f (\tfrac{\jb}{n}) 
         =     6 \sum_{\jb\in \HH_n } 
                   \left (1 - c_\jb^{(n)} \right) f (\tfrac{\jb}{n}). 
\end{align*}
Consequently, we have 
\begin{align*}
 \sum_{\jb \in \HH^*_n}  c_\jb^{(n)} f(\tfrac{\jb}{n}) = &  
    \sum_{\jb \in \HH_n}  c_\jb^{(n)} f(\tfrac{\jb}{n}) +
    \sum_{\jb \in \HH^*_n\setminus \HH_{n}}  c_\jb^{(n)} f(\tfrac{\jb}{n}) \\
=&   \sum_{\jb \in \HH_n}  c_\jb^{(n)} f(\tfrac{\jb}{n}) +
   \sum_{ \jb \in \HH_n  }   \left(1-c_\jb^{(n)}\right) f(\tfrac{\jb}{n})
      =  \sum_{\jb \in \HH_n}  f(\tfrac{\jb}{n}). 
\end{align*} 
Replacing $f$ by $f \bar g$, we have proved that $\langle f, g \rangle_n = 
\langle f, g \rangle_n^*$ for $H$-periodic functions, which clearly applies to
functions in $\CH_n$. 
\end{proof}


\subsection{Interpolation on the hexagon} 

First we note that by \eqref{kAx},  Theorem \ref{prop:interpolation}
when restricted to the hexagon domain becomes the following: 

\begin{prop} \label{prop:1stH-interpo}
For $n \ge 0$, define 
$$
\CI_n f (\tb): = \sum_{\jb \in \HH_n}f(\tfrac{\jb}{n}) \Phi_n(\tb - \tfrac{\jb}{n}), 
\quad \hbox{where}\quad
   \Phi_n(\tb) = \frac{1}{3n^2} \sum_{\jb \in \HH_n} \phi_\jb(\tb), 
$$
for $f\in C(\Omega)$. Then $\CI_n f \in \CH_n$ and 
$$
   \CI_n f (\tfrac{\jb}{n}) =f (\tfrac{\jb}{n}), \qquad \forall \jb \in \HH_n.
$$
\end{prop}

Furthermore, and perhaps much more interesting, is another interpolation 
operator that works for almost all points in $\HH_n^*$. First we need a few 
more definitions.  We denote by $\partial \HH_n^*:=\HH_n^* \setminus 
\HH_n^\circ$,  the set of points in $\HH_n^*$ that are on the boundary of 
the hexagon $\{\tb: -n \le t_1, t_2, t_3 \le n\}$. 
We further divide $\partial \HH_n^*$ as $\HH_n^\ve \cup \HH_n^\e$, where
$\HH_n^\ve$ consists of the six vertices of $\partial \HH_n$ and 
$\HH_n^\e$ consists of the other points in $\partial \HH_n$, i.e., the edge
points. For $\jb \in \partial \HH_n^*$, we define
$$
 \S_\jb: = \left\{\kb \in \HH_n^*: \tfrac{\kb}{n} \equiv \tfrac{\jb}{n} \mod 3\right\}. 
$$
Because of the tiling property of the hexagon, if $\jb \in \HH_n^\e$ then 
$\S_{\jb}$ contains two points, $\jb$ and $\jb^*\in \HH_n^\e$, where $j^*$ is
on the opposite edge, relative to $j$, of the hexagon; while if $\jb \in \HH_n^\ve$, 
then $\S_\jb$ contains three vertices, $\jb$ and its rotations by the angles 
$2\pi /3$ and $4\pi/3$. 

\begin{thm} \label{prop:H-interpo}
For $n \ge 0$ and $f \in C(\Omega)$, define 
$$
\CI_n^* f (\tb): = \sum_{\jb \in \HH_n^*}f(\tfrac{\jb}{n}) \ell_{\jb,n} (\tb) 
$$
where 
\begin{equation} \label{H-interp-ell}
 \ell_{\jb,n}(\tb) = \Phi_n^*(\tb - \tfrac{\jb}{n})  \quad \hbox{and}\quad
   \Phi_n^*(\tb) = \frac{1}{3n^2} \sum_{\jb \in \HH_n^*} c_\jb^{(n)} \phi_\jb(\tb).
\end{equation}
Then $\CI_n^* f \in \CH_n^*:= \{\phi_\jb: \jb \in \HH_n^*\}$ and 
\begin{equation}\label{H-interpo}
 \CI_n^* f (\tfrac{\jb}{n}) = \begin{cases} f(\tfrac{\jb}{n}), &  
       \jb \in \HH_n^\circ, \\
       \sum_{\kb \in \S_\jb} f(\frac{\kb}{n}), &
       \jb \in \partial \HH_n^* . \end{cases}
\end{equation}
Furthermore, $\Phi_n^*(\tb)$ is a real function and it has a compact formula
\begin{align} \label{Phi_n}
    \Phi_n^*(\tb) = & \frac{1}{3n^2} \left[ \frac{1}{2} \big(\Theta_n(\tb) -
       \Theta_{n-2}(\tb)\big)  \right. \\
     & \qquad
      \left.  - \frac{1}{3} \big( \cos \tfrac{2n\pi}{3} (t_1-t_2)+ \cos \tfrac{2n\pi}{3}(t_2-t_3)
                + \cos  \tfrac{2n\pi}{3} (t_3-t_1)\big) \right]
    \notag
\end{align}
where 
\begin{equation} \label{Theta_n}
 \Theta_n(\tb) : =  \frac{\sin \frac{(n+1)\pi(t_1-t_2)}{3}\sin\frac{(n+1)\pi(t_2-t_3)}{3}
\sin \frac{(n+1)\pi(t_3-t_1)}{3}}{ \sin\frac{\pi(t_1-t_2)}{3}\sin\frac{\pi(t_2-t_3)}{3}
\sin \frac{\pi(t_3-t_1)}{3}}. 
\end{equation}
\end{thm}

\begin{proof} 
First let  $\kb \in \HH_n$.  From the definition of the inner product 
$\la \cdot,\cdot\ra_n^*$, it follows that if $\jb \in \HH_n^\circ$, then 
$
  \ell_{\jb,n}(\tfrac{\kb}{n}) = \la \phi_\kb,\phi_\jb \ra_n^* =\delta_{\kb,\jb}. 
$
If $\jb \in \partial \HH_n^*$, then one of the points in $\S_\jb$ will be in 
$\HH_n$, call it $\jb^*$, then $\ell_{\jb,n}(\tfrac{\kb}{n}) = \delta_{\kb,\jb^*}$. 
As the role of $\jb$ and $\kb$ is symmetric in $\ell_{\jb,n} (\frac{\kb}{n})$,
this covers all cases. Putting these cases together, we have proved that 
\begin{equation} \label{H-ell}
\ell_{\jb,n}(\tfrac{\kb}{n}) = \Phi_n^*(\tfrac{\kb}{n}-\tfrac{\jb}{n}) 
       = \begin{cases} 1 & \hbox{if $\tfrac{\kb}{n} \equiv \tfrac{\jb}{n} \mod 3$} \\
                 0 & \hbox{otherwise} \end{cases}, \qquad \kb,\jb \in \HH_n^*, 
\end{equation}
from which the interpolation \eqref{H-interpo} follows immediately. 
 
By the definition of $c_\jb^{(m)}$, it is easy to see that 
$$
   \Phi_n^*(\tb) = \frac{1}{3n^2} \left[ \frac{1}{2} \big( D_n(\tb) + D_{n-1}(\tb)\big)
    - \frac{1}{6} \sum_{\jb\in \HH_n^\ve} \phi_\jb(\tb) \right],
$$
where $D_n$ is the analogue of the Dirichlet kernel defined by
\begin{equation}\label{Gamma}
  D_n(\tb) : =  \sum_{\jb \in \HH_n^*} \phi_\jb(\tb). 
\end{equation}
By the definitions of $\phi_\kb$ and $\HH_n^\ve$, we have 
$$
  \sum_{\jb\in \HH_n^\ve} \phi_\jb(\tb) = 2 \left[ \cos  \tfrac{2n\pi}{3} (t_1-t_2)
       + \cos  \tfrac{2n\pi}{3} (t_2-t_3)+     \cos  \tfrac{2n\pi}{3} (t_3-t_1)\right],
$$
so that we only have to show that $D_n=\Theta_n- \Theta_{n-1}$, 
which in fact has already appeared in \cite{Sun}. Since an intermediate step 
is needed later, we outline the proof of the identity below. The hexagon domain 
can be partitioned into  three parallelograms as shown in Figure 4. 

\begin{figure}[h]
\centerline{\includegraphics[width=0.48\textwidth]{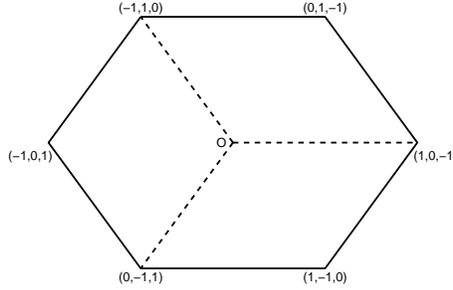}}
\caption{ The hexagon partitioned into three parallelograms.}
 \end{figure}


\noindent
Accordingly, we split the sum over $\HH_n^*$ as 
\begin{equation} \label{compact1}
   D_n(\tb)=  
       1+ \sum_{j_2=-n}^{-1}\sum_{j_1=0}^n \phi_\jb(\tb)+ 
            \sum_{j_3=-n}^{-1}\sum_{j_2=0}^n \phi_\jb(\tb) + 
           \sum_{j_1=-n}^{-1}\sum_{j_3=0}^n \phi_\jb(\tb).
\end{equation}
Using the fact that $t_1+t_2+t_3 =0$ and $j_1+j_2+j_3=0$,  each of the three 
sums can be summed up explicitly. For example, 
\begin{align} \label{compact2}
&  \sum_{j_2=-n}^{-1}\sum_{j_1=0}^n \phi_\jb(\tb) \\
&  \quad 
 =  \frac{\e^{- \frac{\pi i}{3} (n+1) (t_3-t_2)} \sin \frac{n\pi}{3}(t_2-t_3)}
          {\sin \frac{\pi}{3}(t_2-t_3)} 
         \frac{\e^{- \frac{\pi i}{3} n (t_1-t_3)} \sin \frac{(n+1)\pi}{3}(t_1-t_3)}
          {\sin \frac{\pi}{3}(t_1-t_3)}.  \notag
\end{align} 
The formula $D_n=\Theta_n- \Theta_{n-1}$ comes from summing 
up these terms and carrying out simplification using the fact that
if $t_1+t_2+t_3 =0$, then 
\begin{align} \label{trig_identity}
\begin{split}
   \sin 2 t_1 + \sin 2 t_2 + \sin 2 t_3 & = -4 \sin t_1 \sin t_2 \sin t_3, \\
   \cos 2 t_1 + \cos 2 t_2 + \cos 2 t_3 & = 4 \cos t_1 \cos t_2 \cos t_3 -1, 
\end{split}
\end{align}
 which can be easily verified. 
\end{proof}  

The compact formula of the interpolation function allows us to estimate the 
operator norm of $\CI_n^*: C(\overline{\Omega}) \mapsto C(\overline{\Omega})$
in the uniform norm, which is often referred to as the Lebesgue constant. 

\begin{thm} \label{I-norm}
Let $\|\CI_n^*\|_\infty$ denote the operator norm of $\CI_n^*: C(\overline{\Omega}) 
\mapsto C(\overline{\Omega})$. Then there is a constant $c$, independent of 
$n$, such that 
$$
     \|\CI_n^* \|_\infty \le c (\log n )^2. 
$$
\end{thm}

\begin{proof}
A standard procedure shows that 
$$
  \|\CI_n^*\|_\infty = \max_{\tb \in \Omega} 
          \sum_{\jb \in \HH_n^*} | \Phi_n^*(\tb - \tfrac{\jb}{n})|.  
$$
By \eqref{Phi_n} and the fact that $D_n = \Theta_n - \Theta_{n-1}$, it is 
sufficient to show that 
$$
     \frac{1}{3 n^2} \max_{\tb \in \Omega}
    \sum_{\jb \in \HH_n^*} | D_n(\tb - \tfrac{\jb}{n})|\le c (\log n )^2, \qquad 
    n \ge 0.
$$ 
By \eqref{compact1} and \eqref{compact2}, the estimate can be further 
reduced to show that 
\begin{align} \label{estimate}
  J_n(\tb) & := \frac{1}{3n^2} \sum_{\jb \in \HH_n^*} \left \vert 
       \frac{ \sin \frac{n\pi}{3}(t_2-t_3 - \tfrac{j_2-j_3}{n})} 
          {\sin \frac{\pi}{3}(t_2-t_3 -\tfrac{j_2-j_3}{n})} 
         \frac{ \sin \frac{(n+1)\pi}{3}(t_1-t_3- \tfrac{j_1-j_3}{n})}
          {\sin \frac{\pi}{3}(t_1-t_3- \tfrac{j_1-j_3}{n})}\right \vert  \\
       & :=  \sum_{\jb \in \HH_n^*} F_n (\tb - \tfrac{\jb}{n}) \notag
\end{align}
and two other similar sums obtained by permuting $t_1,t_2,t_3$ are 
bounded by  $c (\log n )^2$ for all $\tb \in \Omega$. The three sums are 
similar, we only work with \eqref{estimate}. 

Fix a $\tb \in \Omega$. If $s_1+s_2+s_3 =0$, then the equations 
$s_2-s_3 = t_2-t_3$ and $s_1-s_3=t_1-t_3$ has a unique solution 
$s_2= t_2$ and $s_3 = t_3$. In particular, there can be at most one
$\jb$ such that the denominator of $F_n(\tb - \tfrac{\jb}{n})$ is zero. 
For $\tb \in \overline{\Omega}$, setting $s_1 = (t_1 - t_3 )/3 = (2t_1+t_2)/3$ and 
$s_2 = (t_2- t_3 )/3=(t_1+2t_2)/3$, it is easy to see that $-1 \le s_1,s_2 \le 1$.
The same consideration also shows that $-3n \le j_1-j_3, j_2 - j_3 \le 3n$ for 
$\jb \in \HH_n^*$. Consequently, enlarging the set over which the summation is
taken, it follows that 
$$
  J_n(\tb) \le  \frac{1}{3n^2}  
     \sum_{k_2=-3 n}^{3n} \left \vert 
       \frac{ \sin n\pi(s_2 - \tfrac{k_2}{3 n})} 
          {\sin \pi (s_2 -\tfrac{k_2}{3 n})} \right\vert
     \sum_{k_1=-3 n}^{3n} \left \vert   \frac{\sin (n+1)\pi(s_1- \tfrac{k_1}{3n})}
          {\sin \pi(s_1- \tfrac{k_1}{3 n})}\right \vert. 
$$
For $s \in [-1,1]$, a well-known procedure in one variable shows that 
$$
\max_{s \in [-1,1]} 
   \frac{1}{n} \sum_{k=-3 n}^{3n} \left \vert  \frac{ \sin n\pi(s - \tfrac{k}{3 n})} 
          {\sin \pi (s -\tfrac{k}{3 n})} \right\vert \le c \log n,
$$
from which the stated results follows immediately. 
\end{proof}


\section{Discrete Fourier analysis on the triangle} 
\setcounter{equation}{0}

Working with functions invariant under the isometries of the hexagon lattice,
we can carry out a Fourier analysis on the equilateral triangle based on the
analysis on the hexagon. 

\subsection{Generalized sine and cosine functions}

The group $\G$ of isometries of the hexagon lattice is generated by the reflections 
in the edges of the equilateral triangles inside it, which is the reflection group 
$\A_2$. In homogeneous coordinates, the three reflections $\sigma_1,
\sigma_2,\sigma_3$ are defined by
$$
  \tb  \sigma_1 :=  -(t_1,t_3,t_2),  \quad \tb \sigma_2 := -(t_2,t_1,t_3),  
     \quad \tb\sigma_3:= -(t_3,t_2,t_1).  
$$
Because of the relations $\sigma_3 = \sigma_1\sigma_2\sigma_1
=\sigma_2\sigma_1\sigma_2$, the group $\G$ is given by 
$$
\G = \{1, \sigma_1,\sigma_2,\sigma_3, \sigma_1\sigma_2,\sigma_2\sigma_1\}.
$$
For a function $f$ in homogeneous coordinates, the action of the group $\G$ 
on $f$ is defined by $\sigma f(\tb) = f (\tb \sigma)$, $\sigma \in \G$. Following 
\cite{K}, we call a function $f$ {\it invariant} under $\G$ if $\sigma f =f$ for all 
$\sigma \in \G$, and call it {\it anti-invariant} under $\G$ if $\sigma f = \rho(\sigma)f$ 
for all $\sigma \in \G$, where $\rho(\sigma) = -1$ if $\sigma = \sigma_1,\sigma_2,\sigma_3$, and $\rho(\sigma)=1$ for other elements in $\G$. The following 
proposition is easy to verify (\cite{K}). 

\begin{prop} 
Define two operators $\CP^+$ and $\CP^-$ acting on $f(\tb)$ by 
\begin{equation} \label{CP^+}
\CP^\pm f(\tb) =  \frac{1}{6} \left[f(\tb) + f(\tb \sigma_1\sigma_2)+f(\tb \sigma_2\sigma_1)
        \pm  f(\tb \sigma_1) \pm  f(\tb \sigma_2)  \pm  f(\tb \sigma_3) \right]. 
\end{equation}
Then the operators $\CP^+$ and $\CP^-$ are projections from the class of
H-periodic functions onto the class of invariant, respectively anti-invariant 
functions.  
\end{prop}

Recall that  $\phi_\kb(\tb)= \e^{\tfrac{2 \pi i}{3} \kb \cdot \tb}$.  Following 
\cite{Sun}, we shall call the functions 
\begin{align*} 
 \TC_\kb(\tb) :=  \CP^+ \phi_\kb(\tb)  = & 
    \frac{1}{6} \left[ \phi_{k_1,k_2,k_3}(\tb)+ \phi_{k_2,k_3,k_1}(\tb)+
     \phi_{k_3,k_1,k_2}(\tb) \right.\\
  & \left.+ \phi_{-k_1,-k_3,-k_2}(\tb) +\phi_{-k_2,-k_1,-k_3}(\tb) 
     +\phi_{-k_3,-k_2,-k_1}(\tb)\right], \\
 \TS_\kb(\tb): = \frac{1}{i} \CP^- \phi_\kb(\tb) = &
    \frac{1}{6i } \left[ \phi_{k_1,k_2,k_3}(\tb)+ \phi_{k_2,k_3,k_1}(\tb)+
     \phi_{k_3,k_1,k_2}(\tb) \right.\\
& \left.- \phi_{-k_1,-k_3,-k_2}(\tb) -\phi_{-k_2,-k_1,-k_3}(\tb) 
      -\phi_{-k_3,-k_2,-k_1}(\tb)\right] 
\end{align*}
a {\it generalized cosine} and a {\it generalized sine}, respectively,  where the 
second equal sign follows from the fact that $\phi_\kb (\tb \sigma) = 
 \phi_{\kb \sigma}(\tb)$ for every $\sigma \in \G$. It follows that
 $\TC_\kb$ and $\TS_\kb$ are invariant and anti-invariant, respectively. 

These two functions appear as eigenfunctions of the Laplacian on the 
triangle and they  have been studied in \cite{K,Mc1,Mc2,P1,P2,Sun,LS} 
under various coordinate systems. In particular, they are shown to share 
many properties of the classical sine and cosine functions in 
homogeneous coordinates in \cite{Sun,LS}. Some of those properties that 
we shall need will be recalled below. The advantage of homogeneous 
coordinates lies in the symmetry of the formulas. For example,  the equations 
$t_1+t_2+t_3=0$ and $k_1+k_2+k_3=0$ can be used to write $\TC_\kb$ and 
$\TS_\kb$ in several useful forms, such as
\begin{align}
\TC_\kb(\tb) = & \frac{1}{3} \left[ \e^{\frac{i\pi}{3}(k_2-k_3)(t_2-t_3)}\cos k_1\pi t_1  \label{TC_cos} \right.\\
     & \left. +    \e^{\frac{i\pi}{3}(k_2-k_3)(t_3-t_1)}\cos k_1\pi t_2  
   +\e^{\frac{i\pi}{3}(k_2-k_3)(t_1-t_2)}\cos k_1\pi t_3\right],  \notag\\ 
\TS_\kb(\tb) = &  \frac{1}{3} \left[
   \e^{\frac{i\pi}{3}(k_2-k_3)(t_2-t_3)}\sin k_1\pi t_1\label{TS_sin} \right.\\
     &  \left. +    \e^{\frac{i\pi}{3}(k_2-k_3)(t_3-t_1)}\sin k_1\pi t_2  
   +\e^{\frac{i\pi}{3}(k_2-k_3)(t_1-t_2)}\sin k_1\pi t_3\right], \notag
\end{align}
and similar formulas derived from the permutations of $t_1,t_2, t_3$. In particular, 
it follows from \eqref{TS_sin} that $\TS_\kb(\tb) \equiv 0$ whenever $\kb$ contains
one zero component. 

For invariant functions, we can use symmetry to translate the results
on the hexagon to one of its six equilateral triangles. We shall choose the
triangle as 
\begin{align} \label{Delta}
   \Delta := & \{(t_1,t_2,t_3) : t_1 + t_2 + t_3 =0,   0 \le t_1,  t_2, -t_3 \le 1\}\\
           = & \{(t_1,t_2): t_1, t_2 \ge 0, \, t_1+t_2 \le 1\}. \notag
\end{align}
The region $\Delta$ and its relative position in the hexagon is depicted 
in Figure 5.

\begin{figure}[h]
\centerline{ \includegraphics[width=0.4\textwidth]{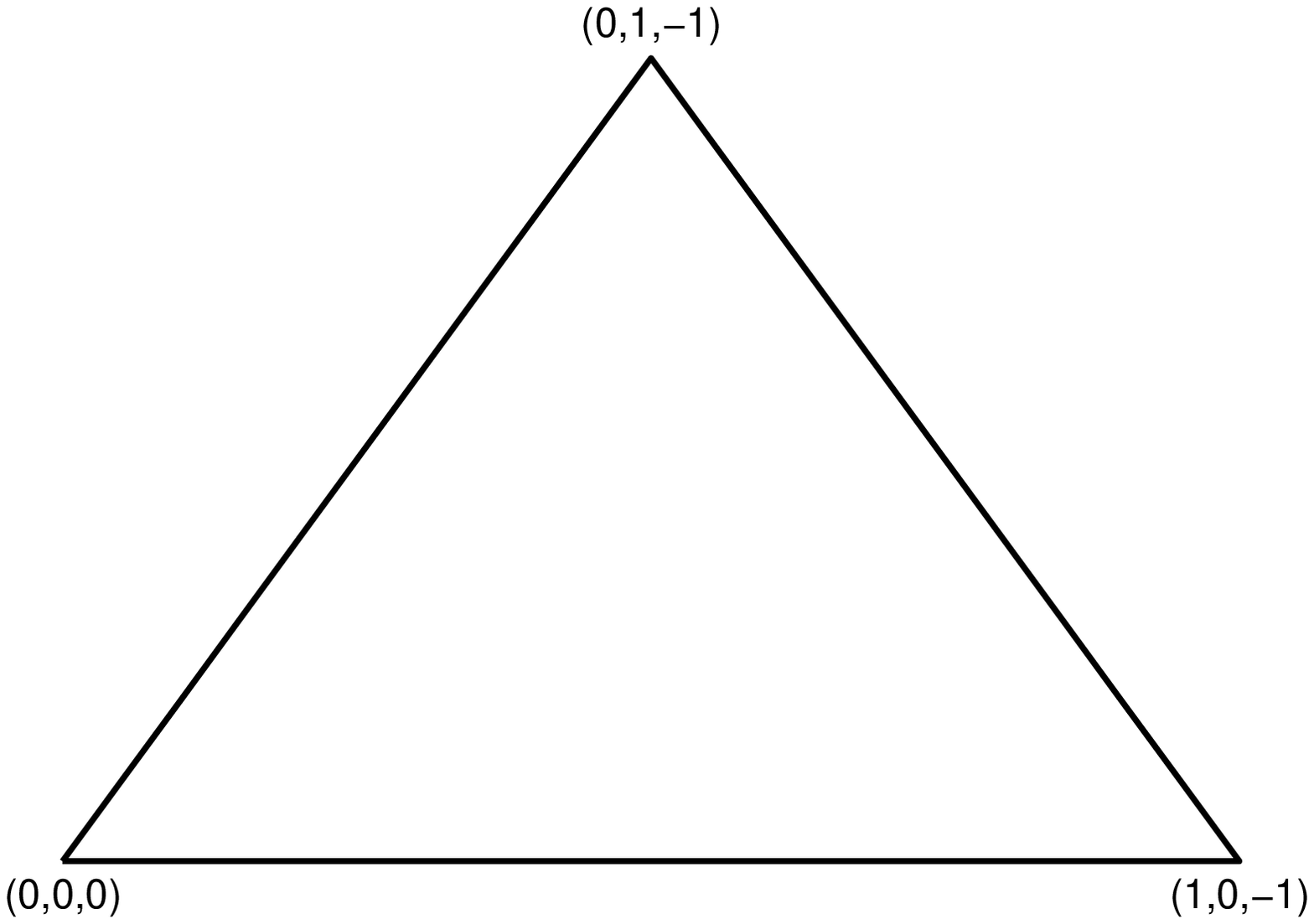} \quad
\includegraphics[width=0.56\textwidth]{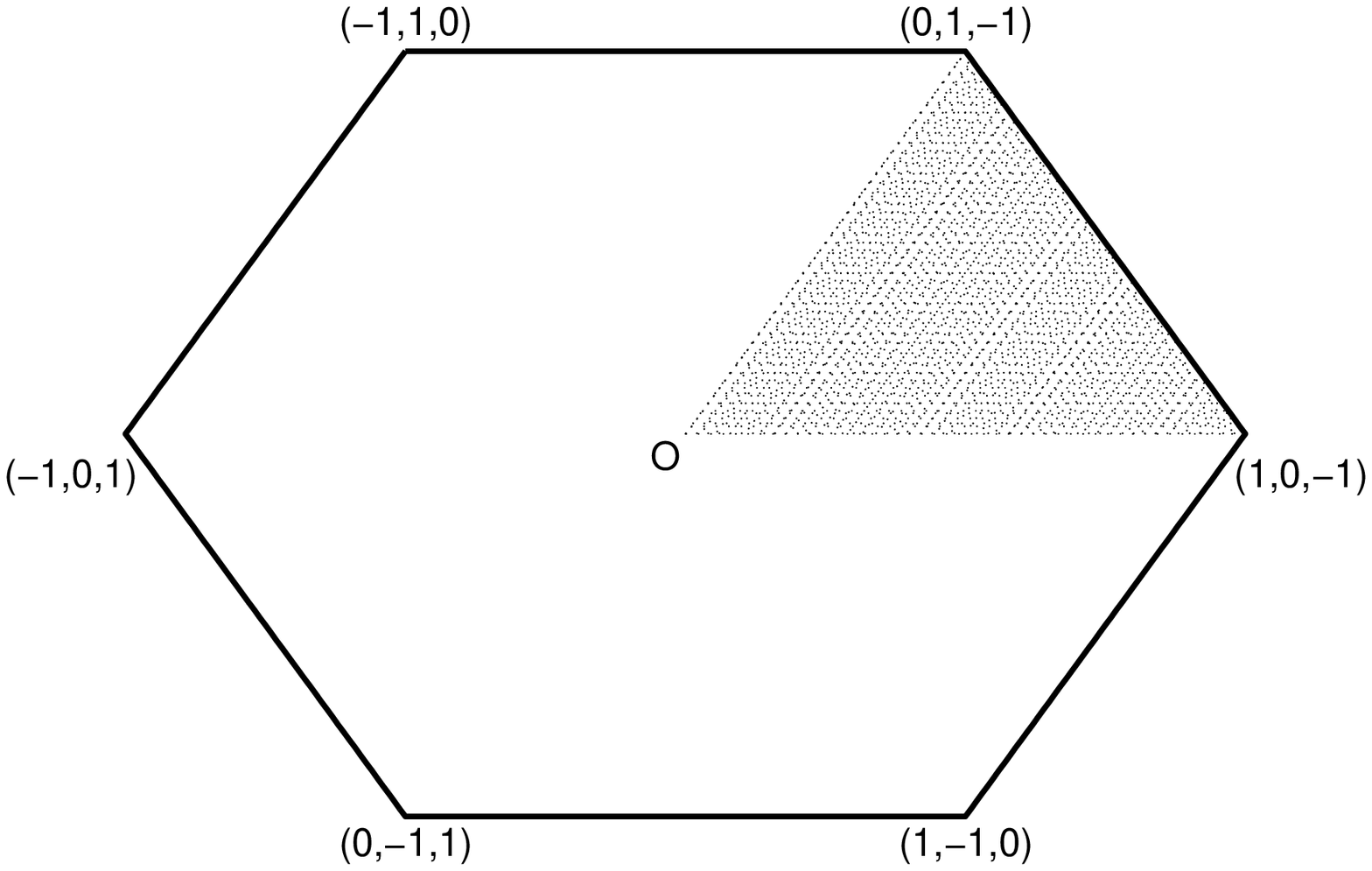} }
\caption{ The equilateral triangle $\Delta$ in the hexagon.}
\end{figure}

\noindent
The inner product on $\Delta$ is defined by 
$$
   \langle f, g \rangle_\Delta :=  
       \frac{1}{|\Delta|}\int_\Delta f(\tb)\overline{g(\tb)} d\tb
      = 2 \int_\Delta f(t_1,t_2) \overline{g(t_1,t_2)} dt_1 dt_2.
$$
If $f \bar g$ is invariant under $\G$, then it is easy to see that $\la f,g\ra_H = 
\la f, g \ra_\Delta$. 

When $\TC_\kb$ are restricted to the triangle $\Delta$, we only need to consider
a subset of $\kb \in \HH$. In fact, since $\TC_{\kb \sigma}(\tb) = \TC_\kb(\tb \sigma) 
=  \TC_\kb(\tb)$ for $\tb \in \Delta$ and $\sigma \in \G$, we can restrict  
$\kb$ to the index set
\begin{equation}\label{Lambda}
  \Lambda: = \{\kb \in \HH: k_1 \ge 0, k_2 \ge 0, k_3 \le 0\}. 
\end{equation}
We define the space $\CT$ as the collection of finite linear combinations of 
elements in  $\{\TC_\kb, \TS_\kb: \kb \in \Lambda\}$; that is, 
\begin{equation}\label{CT-space}
  \CT = \left\{\sum_{\kb \in \Xi} (a_\kb \TC_\kb + b_\kb \TS_\kb)
  : a_\kb, b_\kb  \in \CC, \Xi \subset \Lambda,
    \# \Xi<\infty \right\}.
\end{equation}
Since $\TS_\kb(\tb) = 0$ whenever $\kb$ has a zero component, 
$\TS_\kb (\tb)$ are defined only for $\kb \in \Lambda^\circ$, where
$$
     \Lambda^\circ: = \{\kb \in \HH: k_1 > 0, k_2 > 0, k_3 < 0\},
$$
which is the set of the interior points of $\Lambda$.  These functions are 
orthogonal on $\Delta$.

\begin{prop}\label{prop:trig-ortho}
For $\kb, \jb \in \Lambda$, 
\begin{equation}\label{TC-ortho}
     \la\TC_\kb, \TC_{\jb} \ra_\Delta  =  \delta_{\kb, \jb} \begin{cases}
       1, & \kb = 0,\\
       \frac{1}{3}, & \hbox{$\kb \in \Lambda\setminus \Lambda^\circ$, $\kb \ne 0$},\\
       \frac{1}{6}, & \hbox{$\kb \in \Lambda^\circ$},  \end{cases} 
\end{equation}
and for $\kb,\jb \in \Lambda^\circ$, 
\begin{equation}\label{TS-ortho}
  \la\TS_\kb, \TS_{\jb} \ra_\Delta  = \frac{1}{6} \,\delta_{\kb, \jb}.
\end{equation}
\end{prop}

This follows easily from $\la f, g \ra_H = \la f, g\ra_\Delta$, which holds if
$f \bar g$ is invariant, the fact that 
$\TC_\kb$ is invariant and $\TS_\kb $ is anti-invariant, and the orthogonality
of $\phi_\kb$ in \eqref{H-expOrth}.  The norms of $\TC_\kb$ are divided into 
three cases, since $\TC_0(\tb) =1$ and
$$
    \TC_\kb(\tb) = \frac{1}{3} \left[\phi_{k_1,k_2,k_3}(\tb)+\phi_{k_2,k_3,k_1}(\tb)
        +\phi_{k_3,k_1,k_2}(\tb) \right], \qquad k_1k_2k_3 =0. 
$$


\subsection{Discrete inner product on the triangle}

Using the fact that $\TC_\kb$ and $\TS_\kb$ are invariant and anti-invariant
under $\G$ and the orthogonality of $\phi_\kb$ with respect to the symmetric inner product $\la \cdot, \cdot\ra_n^*$, we can deduce a discrete orthogonality for the 
generalized cosine and sine functions. For this purpose, we define 
$$
  \Lambda_n := \{\kb \in \Lambda:   - k_3 \le n\}
      = \{(k_1,k_2) \in \ZZ^2: k_1 \ge 0, k_2 \ge 0, k_1+k_2 \le n\},
$$
and denote by $\Lambda^\circ, \Lambda_n^\e, \Lambda_n^\ve$
the set of points in the interior of $\Lambda_n$, on the boundary of 
$\Lambda_n$ but not vertices, and at vertices, respectively; that is,
\begin{align*}
 \Lambda_n^\circ & := \{\kb \in \Lambda_n:  k_1>0, k_2>0, - k_3 < n\}, \\
 \Lambda_n^\ve  & := \{(n,0,-n), (0,n,-n),(0,0,0)\}, \\
 \Lambda_n^\e & := \{\kb \in \Lambda_n: \kb \notin \Lambda_n^\circ 
   \cup\Lambda^\ve\}.
\end{align*}
We then define by $\CTC_n$ and $\CTS_n$ the subspace of $\CT$ defined by 
$$
     \CTC_n  = \sspan\{ \TC_\kb:  \kb \in \Lambda_n\} \quad \hbox{and}\quad
        \CTS_n  = \sspan\{ \TS_\kb:  \kb \in \Lambda_n^\circ\},
$$
respectively. 

\begin{thm} \label{thm:d-inner-triangle}
Let the discrete inner product $\la \cdot,\cdot \ra_{\Delta,n}$ be defined by
$$
  \la f,g \ra_{\Delta,n} = \frac{1}{3n^2} \sum_{\jb \in \Lambda_n} 
        \lambda_{\jb}^{(n)} f(\tfrac{\jb}{n}) \overline{g(\tfrac{\jb}{n})},
$$
where 
\begin{equation}\label{lambda_j}
 \lambda_\jb^{(n)} : =  \begin{cases}
       6, & \jb \in \Lambda^\circ,\\
       3 , & \jb \in \Lambda^\e,\\
       1, & \jb \in \Lambda^\ve. 
       \end{cases} 
\end{equation} 
Then 
\begin{equation} \label{d-inner_triangle}
    \la f, g \ra_{\Delta} = \la f,g\ra_{\Delta, n}, \qquad  \hbox{ $f, g \in \CTC_n$}.       
\end{equation}
\end{thm}

\begin{proof}
We shall deduce the result from the discrete inner product 
$\la \cdot, \cdot \ra_n^*$ defined in \eqref{Symm-inner}. Let $\jb \G$ denote 
the orbit of $\jb$ under $\G$, which is the set $\{\jb \sigma: \sigma \in \G\}$. 
It is easy to see that for $\jb \in \HH$, 
\begin{align} \label{jbG}
     |\jb\G|= \begin{cases} 6 &\,  \hbox{if $j_1j_2j_3 \ne 0$} \\
                 3 & \, \hbox{if $j_1j_2j_3 = 0$ but $\jb \ne (0,0,0)$ } \\
                 1 & \hbox{if $\jb = (0,0,0)$}  \end{cases}. 
\end{align}
Assume that $f$ is invariant. Using \eqref{jbG}, we see that 
\begin{align*}
\sum_{\jb \in \HH_n^*, j_1j_2j_3\ne 0} c_\jb^{(n)} f(\tfrac{\jb}{n}) & = 
  6 \sum_{\jb \in \Lambda_n^\circ} c_{\jb}^{(n)} f(\tfrac{\jb}{n}) 
 + 6 \sum_{\jb \in \Lambda_n^\e, j_1j_2j_3\ne 0} 
     c_{\jb}^{(n)} f(\tfrac{\jb}{n}) \\
   & =   \sum_{\jb \in \Lambda_n, j_1j_2j_3 \ne 0} 
    \lambda_\jb^{(n)} f(\tfrac{\jb}{n}),
\end{align*}
and 
\begin{align*}
\sum_{\jb \in \HH_n^*,j_1j_2j_3 =0} c_\jb^{(n)} f(\tfrac{\jb}{n})
  & = f(0,0,0) + 3 \sum_{j\in \Lambda^\e, j_1j_2j_3=0} c_{\jb}^{(n)} f(\tfrac{\jb}{n})      
     + 3 \sum_{j\in \Lambda^\ve, \jb\ne 0} c_{\jb}^{(n)} f(\tfrac{\jb}{n})  \\    
& = \sum_{\jb \in \Lambda_n,j_1j_2j_3 =0}  \lambda_\jb^{(n)}  f(\tfrac{\jb}{n}). 
\end{align*}
Replacing $f$ by $f\bar g$ and adding these two sums, we have proved
that $\la f,g\ra_n^* = \la f, g \ra_{\Delta, n}$ for invariant $f \bar g$, which
applies to $f, g \in \CTC_{n-1}$. Thus, the stated result follows immediately from 
Theorem \ref{prop: Symm-inner} and the fact that 
$\la f,g\ra_H = \la f, g \ra_{\Delta}$. 
\end{proof}

The proof of the above proposition also applies to $f, g \in \CTS_n$ as 
$f \bar g$ is invariant if both $f$ and $g$ are anti-invariant. Noticing also 
that $\CTS_\jb(\tfrac{\jb}{n}) = 0$ when $\jb \in \Lambda^\e$, and hence 
we conclude the following result. 

\begin{prop}\label{thm:d-inner-triangle2}
Let the discrete inner product $\la \cdot,\cdot \ra_{\Delta^\circ,n}$ be 
defined by
$$
  \la f,g \ra_{\Delta^\circ,n} = \frac{2}{n^2}
     \sum_{\jb \in \Lambda_n^\circ} 
        f(\tfrac{\jb}{n}) \overline{g(\tfrac{\jb}{n})}. 
$$
Then 
\begin{equation} \label{d-inner_triangle2}
    \la f, g \ra_{\Delta} = \la f,g\ra_{\Delta^\circ, n}, \qquad  f, g \in \CTS_n.       
\end{equation}
\end{prop}

\begin{figure}[h]
\centerline{ \includegraphics[width=0.61\textwidth]{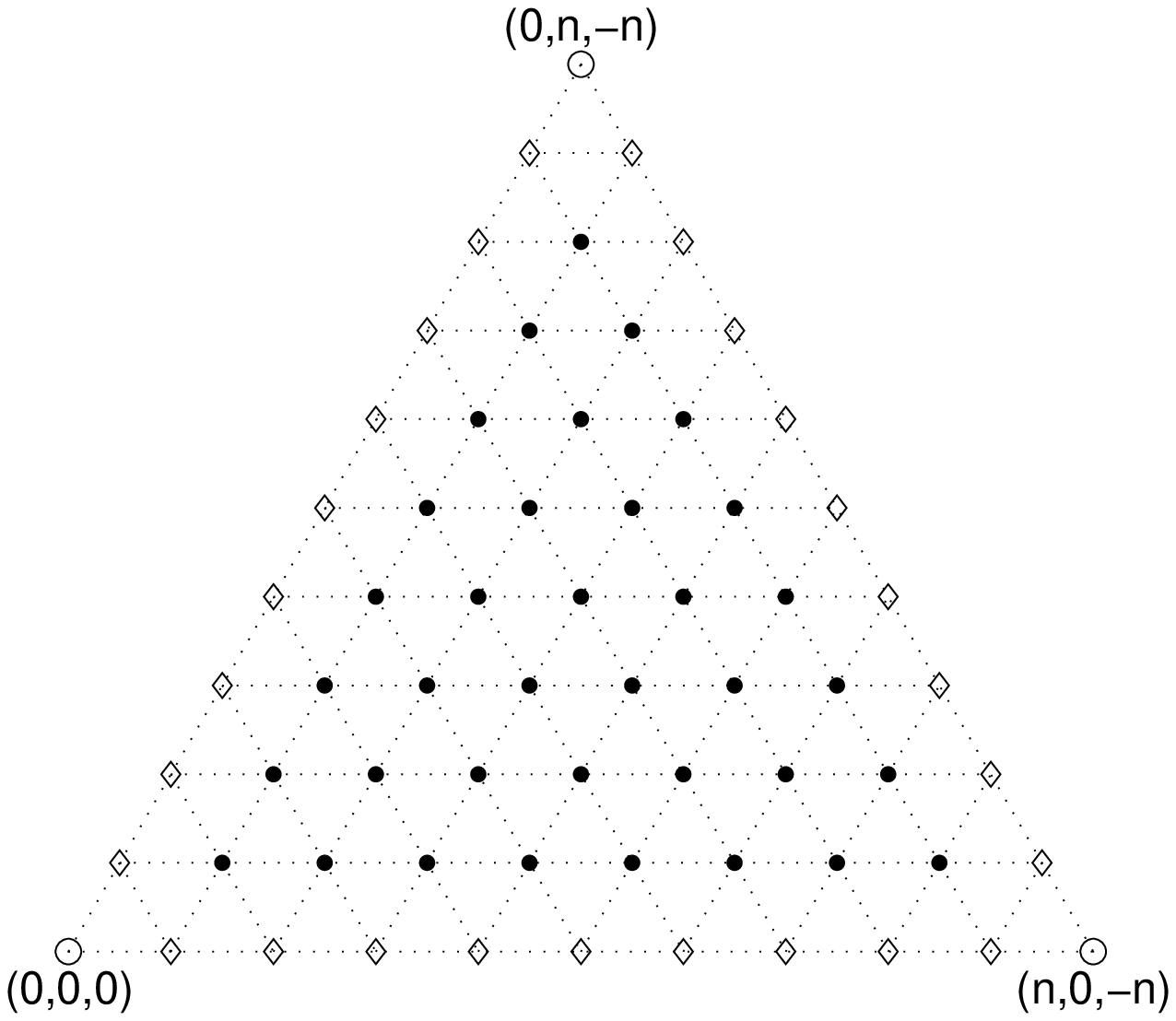} \,
\includegraphics[width=0.48\textwidth]{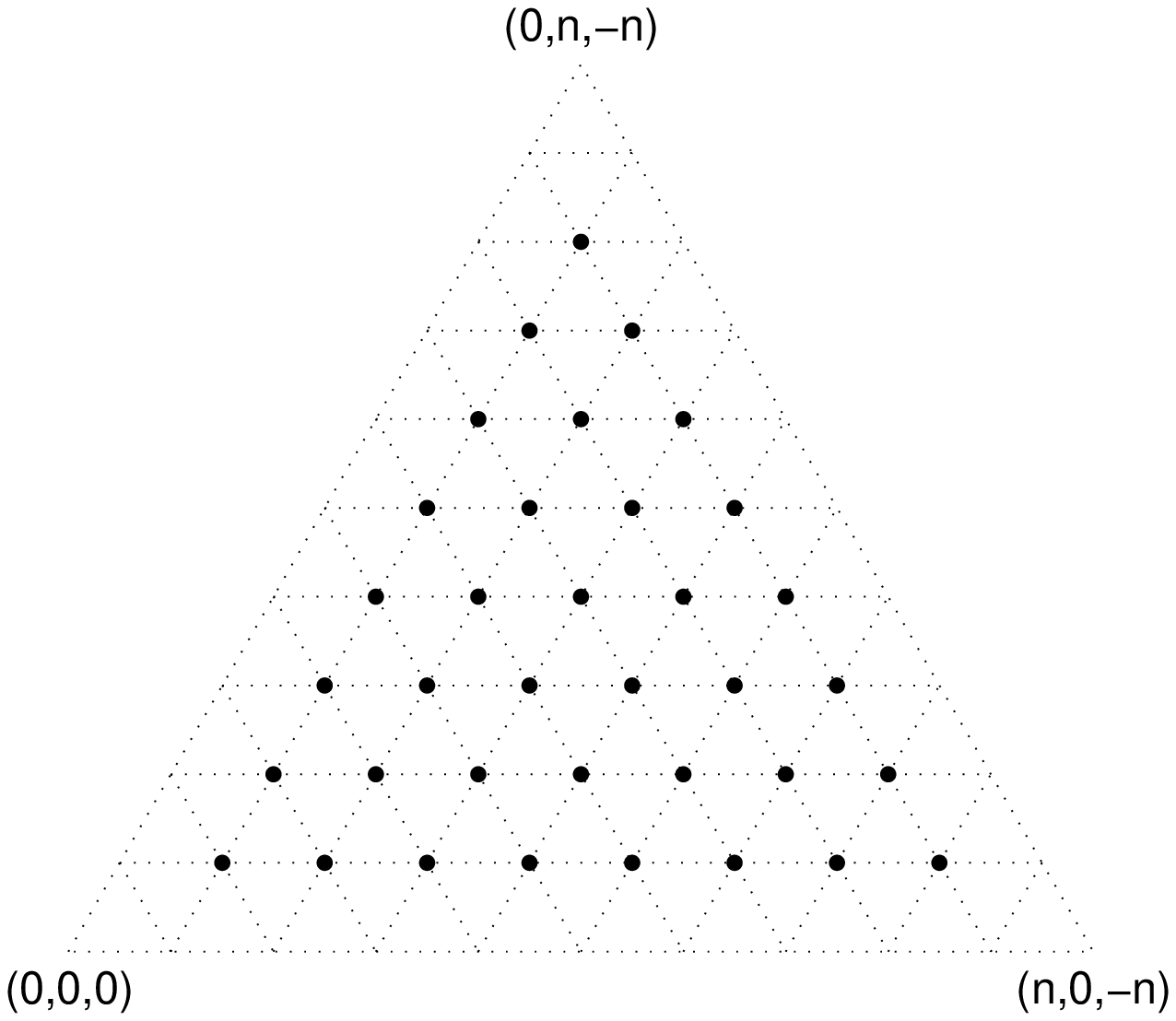} }
\caption{ Left: $\Lambda_n$ in which $\bullet$, $\diamond$ 
and $\circ$ correspond to points of different coefficients in 
$\la\cdot, \cdot\ra_n^*$.  Right:  $\bullet$ indicates points 
in $\Lambda^\circ$.}
 \end{figure}

The relation \eqref{d-inner_triangle} can be regarded as a cubature formula
$$
 \frac{1}{|\Delta|} \int_{\Delta} f(\tb) d\tb = \frac{1}{3n^2}  \sum_{\jb \in \Lambda_n} 
      \lambda_\jb^{(n)}f(\tfrac{\jb}{n})
$$
which is exact for $f = g \bar h$, $g, h \in \CTC_n$. It turns out, however, that 
much more is true. We state the following result in terms of the 
triangle in $(t_1,t_2)$.

\begin{thm} \label{thm:trigCuba}
For $n \ge 0$, the cubature formula 
\begin{equation}\label{trigCuba}
  2  \int_{\Delta} f(t_1,t_2) dt_1 dt_2 = \frac{1}{3n^2}
      \sum_{j_1=0}^n\sum_{j_2=0}^{j_1} 
          \lambda_\jb^{(n)} f(\frac{j_1}{n},\frac{j_2}{n})
\end{equation} 
is exact for all $f \in \CTC_{2n-1}$.  
\end{thm}

In order to establish such a result, we will need to study the structure of 
$\CTC_n$, which, and the proof of the theorem, will be given in 
Section \ref{Sec:Chebyshev}. One may note that the formulation of the result
resembles a Gaussian quadrature. There is indeed a connection which will also 
be explored in Section \ref{Sec:Chebyshev}. 


\subsection{Interpolation on the triangle}
One way to deduce the result on interpolation is by making use of the
orthogonality of generalized trigonometric functions with respect to the 
discrete inner product, as shown in the first theorem below.  Recall the 
operator $\CP^\pm$ defined in \eqref{CP^+}. 

\begin{thm} \label{prop:1st-interpo-Dela}
For $n\ge 0$ and $f\in C(\Delta)$ define 
$$
 \CL_n f (\tb) : = \sum_{\jb \in \Lambda_n^\circ} f(\tfrac{\jb}{n}) 
      \ell_{\jb,n}^{\circ} (\tb), \qquad 
       \ell_{\jb,n}^\circ (\tb) = \frac{2}{n^2} \sum_{\kb \in \Lambda_n^\circ} \TS_{\kb}(\tb) 
         \overline{ \TS_{\kb}(\tfrac{\jb}{n}) }. 
$$
Then $\CL_n$ is the unique function in $\CTS_n$ that satisfies 
$$
  \CL_nf(\tfrac{\jb}{n}) = f(\tfrac{\jb}{n}), \qquad \jb \in \Lambda_n^\circ. 
$$
Furthermore, the fundamental interpolation function $\ell_{\jb,n}^\circ$ is 
real and satisfies 
\begin{align*}
 \ell_{\jb,n}^\circ(\tb) = \frac{1}{3 n^2} 
           \CP^-_\tb \left[ \Theta_{n-1}(\tb - \tfrac{\jb}{n}) -
                  \Theta_{n-2}(\tb - \tfrac{\jb}{n}) \right], 
\end{align*}
where $\CP_\tb^-$ means that the operator $\CP^-$ is acting on the the 
variable $\tb$ and  $\Theta_n$ is defined in \eqref{Theta_n}. 
\end{thm} 

\begin{proof}
By \eqref{TS-ortho} and \eqref{d-inner_triangle2}, 
$\la \TS_\jb, \TS_\kb \ra_{\Delta^\circ,n} = \tfrac{1}{6} \delta_{\kb,\jb}$,
which shows that $\ell_{\jb,n}^\circ(\tfrac{\kb}{n}) = \delta_{\kb,\jb}$ and
verifies the interpolation condition. Furthermore, since $\TS_{\kb}(\tfrac{\jb}{n})
=\TS_{\jb}(\tfrac{\kb}{n})= 0$ whenever $\kb \in \Lambda_n \setminus 
\Lambda_n^\circ$, and $\TS_\jb(\tb) = 0$ whenever $j_1j_2j_3 =0$, 
it follows from the  definition of $\TS_\jb$ that  
$$
  \ell_{\jb,n}^\circ (\tb)
       = \frac{1}{3n^2} \CP^-_\tb \CP^-_\jb \sum_{\kb \in \HH_n^\circ} \phi_{\kb}(\tb) 
                  \overline{ \phi_{\kb}(\tfrac{\jb}{n})} 
           =  \frac{1}{3n^2} \CP^-_\tb \CP^-_\jb D_{n-1}(\tb - \tfrac{\jb}{n}),
$$
where $D_n$ is defined in \eqref{Gamma}. To complete the proof we use
the fact that if $F$ is an invariant function, then  
$$
   \CP_\tb^- \CP_\sb^- F(\tb - \sb) = \CP_\tb^- F(\tb - s),
$$
which can be easily verified by using the fact that the elements of the 
group $\G$ satisfy $\sigma_i^2 =1$, $i = 1,2,3$ and $\sigma_3  = 
\sigma_1\sigma_2\sigma_1=\sigma_2\sigma_1\sigma_2$.  
\end{proof}

We can derive the result of interpolation on the $\Lambda_n$ using the
same approach, but it is more illustrative to derive it from the interpolation 
on the hexagon given in Theorem \ref{prop:H-interpo}.  

\begin{thm}\label{prop:interpo-Dela}
For $n \ge 0$ and $f\in C(\Delta)$ define 
\begin{equation}\label{interpo-Dela}
  \CL_n^* f(\tb) := 
          \sum_{\jb \in \Lambda_n} f (\tfrac{\jb}{n}) \ell_{\jb,n}^\Delta(\tb), 
     \qquad  \ell_{\jb,n}^\Delta (\tb) = \lambda_{\jb}^{(n)} \CP^+ \ell_{\jb,n}(\tb),
\end{equation}
where $\lambda_\jb^{(n)}$ and $\ell_{\jb,n}$ are defined in \eqref{lambda_j}
and \eqref{H-interp-ell}, respectively. Then $\CL_n^*$ is the 
unique function in $\CTC_n$ that satisfies 
$$
  \CL_n^* f(\tfrac{\jb}{n}) = f(\tfrac{\jb}{n}), \qquad \jb \in \Lambda_n. 
$$
Furthermore, the fundamental interpolation function $\ell_{\jb,n}$ is given by
\begin{align*}
 \ell_{\jb,n}^\Delta(\tb) = \frac{\lambda_\jb^{(n)}}{3n^2} 
       \left[ \frac{1}{2} \CP^+  \big(\Theta_n(\tb - \tfrac{\jb}{n}) 
           -\Theta_{n-2}(\tb - \tfrac{\jb}{n}) \big) - 
   \Re\left\{ \TC_{n,0,-n}(\tb) \overline{\TC_{n,0,-n}(\tfrac{\jb}{n})}  \right\} \right].
\end{align*}
\end{thm}

\begin{proof}
Throughout this proof, write $\ell_\jb = \ell_{\jb,n}$.
Recall the equation \eqref{H-ell}, which states that $\ell_{\jb} (\tfrac{\kb}{n}) = 1$ 
if $\kb \equiv \jb \mod 3$ and $\ell_{\jb} (\tfrac{\kb}{n}) = 0$ otherwise. Let 
$\kb \in \Lambda$. If $\jb \in \Lambda^\circ$, then we have 
$\CP^+ \ell_\jb(\tfrac{\kb}{n}) = \frac{1}{6} \ell_\jb (\tfrac{\kb}{n}) =
[\lambda_\jb^{(n)}]^{-1}\delta_{\kb,\jb}$. If $\jb \in \Lambda^\e$ and $j_1j_2j_3\ne 0$, 
then $\CP^+ \ell_\jb(\tfrac{\kb}{n}) = \frac{1}{6} [\ell_\jb (\tfrac{\kb}{n}) +
\ell_{\jb} (\tfrac{\kb^*}{n})] = \frac{1}{3}  \delta_{\kb,\jb} = 
[\lambda_\jb^{(n)}]^{-1} \delta_{\kb,\jb}$, where $k^*$ is defined as in the
proof of Theorem \ref{prop:H-interpo}. If $\jb \in \Lambda^\e$ and 
$j_1j_2j_3 = 0$, then $\CP^+ \ell_\jb(\tfrac{\kb}{n}) = \frac{1}{3} 
\ell_{\jb} (\tfrac{\kb}{n}) = [\lambda_\jb^{(n)}]^{-1} \delta_{\kb,\jb}$ since 
$\jb \in \HH_n^\circ$. If $\jb \in \Lambda^\ve$ and $\jb \ne (0,0,0)$, then
$\CP^+ \ell_\jb(\tfrac{\kb}{n}) = \frac{1}{3} \sum_{\lb \in \S_\jb} 
\ell_\lb (\tfrac{\kb}{n}) = \delta_{\jb,\kb} =  [\lambda_\jb^{(n)}]^{-1}\delta_{\kb,\jb}$,
where $\S_\jb$ is defined as in the proof of Theorem \ref{prop:H-interpo}. 
Finally, if $\jb = (0,0,0)$, $\CP^+\ell_\jb(\tfrac{\kb}{n}) = \delta_{\kb,\jb}
=[\lambda_\jb^{(n)}]^{-1}\delta_{\kb,\jb}$. Putting these together, we have 
proved that $\ell^\Delta_\jb(\tfrac{\kb}{n}) = \delta_{\kb,\jb}$, which verifies
the interpolation condition. 

The explicit formula of $\ell^\Delta_{\jb}$ follows from \eqref{Phi_n} and 
the fact that 
\begin{align*}
   \CP^+ \sum_{\kb \in \HH_n^\ve} \phi_\kb (\tb- \tfrac{\jb}{n})&  =
      \CP^+ \sum_{\kb \in \HH_n^\ve} \phi_\kb (\tb)\overline{\phi_\kb (\tfrac{\jb}{n})} =
         \sum_{\kb \in \HH_n^\ve} \TC_\kb(\tb) \overline{\phi_\kb (\tfrac{\jb}{n})} \\
      & = 3 \TC_{n,0,-n}(\tb)  \overline{\TC_{n,0,-n} (\tfrac{\jb}{n})} +
             3 \TC_{0,n,-n}(\tb)  \overline{\TC_{0,n,-n} (\tfrac{\jb}{n})},
\end{align*}
as well as the fact that  $\overline{\TC_{0,n,-n}(\tb)} = \TC_{n,0,-n}(\tb)$. 
\end{proof} 

As a consequence of Theorem \ref{I-norm}, we immediately have the following: 

\begin{thm}
Let $\|\CL_n\|$ and $\|\CL_n^*\|$ denote the operator norm of 
$\CL_n : C(\Delta) \mapsto C(\Delta)$ and $\CL_n^* : C(\Delta) \mapsto C(\Delta)$,
respectively. Then there is a constant $c$, independent of $n$, such that
$$
  \| \CL_n \| \le c  (\log n)^2 \quad \hbox{and} \quad  \| \CL_n^* \| \le c (\log n)^2.
$$
In particular, if $f \in C^1(\Omega)$, the class of functions having continuous 
derivatives on $\Omega$, then both $\CL_n f$ and $\CL_n^* f$ converge 
uniformly to $f$ on $\Omega$. 
\end{thm} 

This theorem shows a sharp contrast between interpolation by trigonometric 
functions on $\Lambda_n$ or  $\Lambda_n^\circ$ and the interpolation by 
algebraic polynomials on $\Lambda_n$. From the result in one variable \cite{R}, 
the Lebesgue constant for algebraic polynomial interpolation grows exponentially 
as $n \to \infty$ and, as a consequence, the convergence does not hold in 
$C^1(\Omega)$.


\section{Generalized Chebyshev polynomials and their zeros}\label{Sec:Chebyshev}
\setcounter{equation}{0}

Just like in the classical case of one variable, we can use the generalized
sine and cosine functions to define analogues of Chebyshev polynomials 
of the first and the second kind, respectively, which are orthogonal polynomials
of the two variables. 

\subsection{Generalized Chebyshev polynomials}
To see the polynomial structure among the generalized trigonometric
functions,  we make  a change of variables. Denote 
\begin{align} \label{z} 
z :=  \TC_{0,1,-1}(\tb)  = \tfrac{1}{3} [ \phi_{0,1,-1}(\tb) + \phi_{1,-1,0}(\tb) 
     + \phi_{-1,0,1}(\tb)],
\end{align}
whose real part and imaginary part, after simplification using \eqref{trig_identity},
are
\begin{align} \label{x-y}
\begin{split}
x & = \Re {z}  = \tfrac{4}{3} \cos\tfrac{\pi}{3}(t_2-t_1) \cos\tfrac{\pi}{3}(t_3-t_2) 
  \cos\tfrac{\pi}{3}(t_1-t_3)-\tfrac{1}{3}, \\
y & = \Im {z} =\tfrac{4}{3} \sin\tfrac{\pi}{3}(t_2-t_1) \sin\tfrac{\pi}{3}(t_3-t_2) 
     \sin\tfrac{\pi}{3}(t_1-t_3),
\end{split}
\end{align}
respectively. If we change variables 
$(t_1, t_2) \mapsto (x, y)$, then the region $\Delta$ is mapped onto the region 
$\Delta^*$ bounded by Steiner's hypocycloid, 
\begin{align} \label{Delta*}
\Delta^* 
   = \{(x,y): -3 (x^2+y^2+1)^2 + 8 (x^3-3x y^2)+4\ge 0\}. 
\end{align}

\begin{figure}[h]
\hfill\includegraphics[width=0.4\textwidth,height=0.4\textwidth]{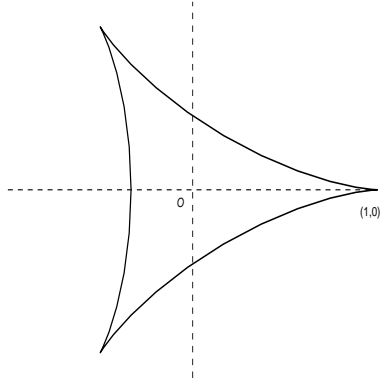}
\hspace*{\fill}
\caption{The region $\triangle^*$ bounded by Steiner's hypocycloid.}
\end{figure}

\begin{defn}
Under the change of variables \eqref{z}, define 
\begin{align*}
T_k^m(z,\bar z) : & = \TC_{k, m-k, - m}(\tb), \qquad 0 \le k \le m, \\
U_k^m(z,\bar z) : & = \frac{\TS_{k+1,m-k+1,-m-2}(\tb)}{\TS_{1,1,-2}(\tb)},
                        \qquad 0 \le k \le m.
\end{align*}
We will also use the notation $T_k^m(\tb)$ and $U_k^m(\tb)$ when it is more 
convenient to work with $\tb$ variable.  
\end{defn}  

These functions are in fact polynomials of degree $m$ in the $z$ and $\bar z$ 
variables since they both satisfy a simple recursive relation.

\begin{prop}
Let $P_k^m$ denote either $T_k^m$ or $U_k^m$. Then 
\begin{equation}\label{conjugT}
  P_{m-k}^m(z, \bar z)  = \overline{P_k^m(z, \bar{z})},   \qquad 0 \le k \le m, 
\end{equation}
and they satisfy the recursion relation
\begin{align} \label{recurT}
     P _k^{m+1} (z,\bar z) = 3 z P_{k}^m(z,\bar z) -
                 P_{k+1}^m(z,\bar z) - P_{k-1}^{m-1}(z,\bar z) 
\end{align}
for $ 0 \le k \le m$ and $m\ge 1$, where we use 
\begin{align*}
&T_{-1}^m(z,\bar z ) = T_1^{m+1}(z,\bar z ), \quad T_{m+1}^m(z,\bar z )
   = T_{m}^{m+1}(z,\bar z ),\\
&U_{-1}^m(z,\bar z ) =  0, \quad  U_m^{m-1}(z,\bar z ) =0.
\end{align*}
In particular, we have 
\begin{align*}
& T_0^0(z,\bar z)=1, \quad T_0^1(z,\bar z)= z,  \quad T_1^1(z,\bar z)= \bar z,\\
& U_0^0(z,\bar z)=1, \quad U_0^1(z,\bar z)= 3z,  \quad U_1^1(z,\bar z)= 3 \bar z.
\end{align*}
\end{prop}

The recursion relation \eqref{recurT} follows from a straightforward computation.
Together, \eqref{recurT} and \eqref{conjugT} determine all $P_k^m$ 
recursively, which show that both $T_k^m$ and $U_k^m$ are polynomials 
of degree $m$ in $z$ and $\bar z$. Furthermore, each family inherits an 
orthogonal relation from the generalized trigonometric functions,  so that they 
are orthogonal polynomials of two variables. They were studied by Koornwinder 
(\cite{K,K2}), who called them the generalized Chebyshev polynomials of the 
first and the second kind, respectively. Moreover, they are special cases of 
orthogonal polynomials with respect to the weighted inner product 
$$
 \la f, g \ra_{w_\alpha}  = c_\alpha \int_{\Delta^*} f(x,y)\overline{g(x,y)}
      w_\alpha(x,y) dxdy,
$$
where $w_\alpha$ is the weight function defined in terms of the Jacobian of
the change of variables \eqref{z},
\begin{align*}
  w_\alpha(x,y) = \left | \frac{\partial(x,y)}{\partial (t_1,t_2)} \right |^{2 \alpha}
  & =  \left[ \frac{16}{27} \pi^2 \sin \pi t_1 \sin \pi t_2 \sin \pi (t_1+ t_2) \right]^{2\alpha}\\
  & =  \frac{4^\alpha}{27^\alpha} \pi^{2\alpha}
      \left[-3(x^2+y^2+1)^2 + 8 (x^3-3x y^2)+4 \right]^\alpha
\end{align*}
and $c_\alpha$ is a normalization constant, 
$c_\alpha := 1/\int_\Delta w_\alpha(x,y) dxdy$. We have, in particular, that
$c_{-\frac{1}{2}}=2$ and $c_{\frac{1}{2}} =81/(8\pi^4)$. Since the change of variables 
\eqref{z} shows immediately that 
$$
 \frac{1}{|\Delta|} \int_{\Delta} f(\tb) d\tb = c_{-\frac{1}{2}}
    \int_{\Delta^*} f(x,y) w_{-\frac12}(x,y) dxdy,
$$
it follows from the orthogonality of $\TC_\jb$ and $\TS_\jb$ that $T_k^m$ and 
$U_k^m$ are orthogonal polynomials with  respect to $w_{-\frac{1}{2}}$ and
$w_{\frac{1}{2}}$, respectively. Furthermore, from Proposition \ref{prop:trig-ortho}
we have
\begin{equation}\label{TCnorm}
   \langle T_k^m, T_j^m\rangle_{ w_{-\frac12}} =d_k^m  \delta_{k,j}, \quad d_k^m :=  
       \begin{cases}  1, & m=k=0, \\
         \frac{1}{3}, & \hbox{$k=0$ or $k=m$, $m >0$},\\
        \frac{1}{6}, & \hbox{$1 \le k \le m-1$,$m >0$,}   \end{cases} 
\end{equation}
and 
\begin{equation}\label{UCnorm}
  \langle U_k^m, U_j^m\rangle_{ w_{\frac12}} = \frac{1}{6}\delta_{k,j}, 
      \qquad 0\le j, k \le m.
\end{equation}

Let $\wt T_k^m = \sqrt{d_k^m} T_k^m$ and $\wt U_k^m = \sqrt{6} U_k^m$.
These are orthonormal polynomials. We introduce the notation 
\begin{equation} \label{PP_n}
         \PP_m = \{P_0^m, P_1^m, \ldots, P_m^m\}
\end{equation}
for either $P_k^m = \wt T_k^m$ or $\wt U_k^m$. The polynomials in $\PP_n$ 
consist of an orthonormal basis for the space $\CV_n$ of orthogonal polynomials 
of degree $m$. Evidently both $T_k^m$ and $U_k^m$ are complex polynomials. 
As a consequence of \eqref{conjugT}, however, it is easy to derive real orthogonal
polynomials with respect to $w_{\pm\frac{1}{2}}(x,y)$ on $\Delta^*$. In fact,
let $\PP_m$ be as in \eqref{PP_n}. Then the polynomials in the set $\PP_m$ 
defined by
\begin{align}\label{realOP}
\begin{split}
  \PP^R_{2m-1}(x,y):= & \{\sqrt{2}\Re\{ P_k^{2m-1}(z,\bar{z})\}, 
     \sqrt{2} \Im\{P_k^{2m-1}(z,\bar{z})\}, 
               0 \le k \le \lfloor m-1\rfloor\} \\
 \PP^R_{2m}(x,y):= & \{P_m^{2m}, \sqrt{2}\Re\{P_k^{2m}(z,\bar{z})\}, 
     \sqrt{2} \Im\{P_k^{2m}(z,\bar{z})\}, 
               0 \le k \le \lfloor m\rfloor\}
\end{split}
\end{align}
are real orthonormal polynomials and form a basis for $\CV_n$ with $n =2m-1$
or $2m$, respectively.


\subsection{Zeros of Chebyshev polynomials and Gaussian cubature formula}
Let $w$ be a nonnegative weight function defined on a compact set 
$\Omega$ in $\RR^2$.  A cubature formula of degree $2n-1$ is a sum of point
evaluations that satisfies 
$$
  \int_\Omega f(x,y) w(x,y) dxdy = \sum_{j=1}^N \lambda_j f(x_j,y_j), 
   \qquad \lambda_j \in \RR
$$
for every $f \in \Pi_{2n-1}^2$, where $\Pi_n^2$ denote the space of polynomials
of degree at most $n$ in two variables. The points $(x_j,y_j)$ are called nodes
of the cubature and the numbers $\lambda_i$ are called weights. It is well-known
that a cubature formula of degree $2n-1$ exists only if $N \ge \dim \Pi_{n-1}^2$.
A cubature that attains such a lower bound is called a Gaussian cubature. 

Let $\PP_n := \{P_k^n: 0 \le k \le n\}$ denote a basis of orthonormal polynomials
of degree $n$ with respect to $w(x,y) dxdy$. It is known that a Gaussian 
cubature exists if and only if $\PP_n$ (every element in $\PP_n$) has 
$\dim \Pi_{n-1}^2$ many common zeros,  which are necessarily real and distinct.
Furthermore, the weight $\lambda_j$  in a Gaussian cubature is given by 
\begin{equation}\label{lambda}
   \lambda_j =   [K_{n-1}(x_j,y_j)]^{-1}, \qquad
           K_n(x,y): = \sum_{k=0}^n\sum_{j=0}^k P_j^k(x)P_j^k(y).
\end{equation}
The definition of $K_n(x,y)$ shows that it is the reproducing kernel of $\Pi_n^2$
in $L^2(\mu)$, so that it is independent of the choice of the basis.  

Unlike the case of one variable, however, Gaussian cubatures do not exist in 
general. We refer to \cite{DX,My,St} for results and further discussions. It is sufficient 
to say that Gaussian cubatures are rare. In fact, the family of weight functions 
identified in \cite{SX} for which the Gaussian cubature exist for all $n$ remains 
the only known example up to now. 

We now test the above theory on the generalized Chebyshev polynomials.
Because of the relation \eqref{realOP}, we can work with the zeros of the complex polynomials in $\PP_n$. Furthermore, it follows from \eqref{conjugT} that the
reproducing kernel 
$$
   K_n(z, w): = \sum_{k=0}^n \PP_k^\tr(z,\bar z) \PP_k(w,\bar w) = 
         \sum_{k=0}^n \sum_{j=0}^k P_j^k(z,\bar z) P_j^k(w,\bar w)
$$
is real and, by \eqref{realOP}, is indeed the reproducing kernel of $\Pi_n^2$.

For the Chebyshev polynomials of the first kind, the polynomials in $\PP_n$ do 
not have common zeros in general.  For example, in the case of $m=2$, we 
have $T_0^2(z,\bar z) = 3z^2-2 \bar z $ and $T_1^2(z,\bar z) =  (3 z \bar z -1)/2$, 
so that the three real orthogonal polynomials of degree $2$ can be given by 
$$
  \left\{ 3(x^2-y^2)-2x, \quad 6 x y +2 y, \quad 3(x^2+y^2)-1 \right\},
$$
which has no common zero at all. Therefore, there is no Gaussian cubature 
for the weigh function $w_{-\frac12}$ on $\Delta^*$.

The Chebyshev polynomials of the second kind, on the other hand, does
have the maximum number of common zeros, so that a Gaussian cubature
exists for $w_{\frac{1}{2}}(\tb)$.  Let 
\begin{align*}
X_n : = \Lambda_{n+2}^\circ 
   =  \{(\tfrac{k_1}{n+2}, \tfrac{k_2}{n+2}) | k_1 \ge 1, k_2 \ge 1, 
        k_1+k_2 \le n+1\}, 
\end{align*}
which is in the interior of $\Delta$, and let $X_n^*$ denote the image of $X_n$ 
in $\Delta^*$ under the mapping \eqref{x-y}. 

\begin{thm}
For the weight function $w_{\frac{1}{2}}$ on $\Delta^*$ a Gaussian cubature 
formula exists; that is, 
\begin{equation} \label{GaussCuba}
  c_{\frac{1}{2}} \int_{\Delta^*} f(x,y) w_{\frac{1}{2}}(x,y) dxdy = 
      \sum_{j_1=1}^{n+1}\sum_{j_2=1}^{n+1-j_1}
         \mu_{j}^{(n)} f(\tfrac{j_1}{n+2}, \tfrac{j_2}{n+2}), \quad \forall f \in \Pi_{2n-1}, 
\end{equation}
where $j = (j_1,j_2)$, 
$$
 \mu_j^{(n)} =  \tfrac{32}{9(n+2)^2} \sin^2 \pi\tfrac{j_1}{n+2}   
    \sin^2 \pi\tfrac{j_2}{n+2} \sin^2 \pi \tfrac{j_1+j_2}{n+2}.    
$$
\end{thm}

\begin{proof} 
The numerator of $U_k^n(\tb)$ is $\TS_{k+1, n-k+1, - n-2} (\tb)$, which is 
zero at all $\frac{j}{n+2}$ according to \eqref{TS_sin}. The first equation 
of \eqref{trig_identity} shows that the denominator of $U_k^n(\tb)$ is 
\begin{equation}\label{TS_112}
    \TS_{1,1,-2}(\tb) = \frac{4}{3} \sin \pi t_1 \sin \pi t_2\sin \pi t_3,
\end{equation}
which vanishes at the integer points on the boundary point  of $\Delta$. 
Consequently, $U_k^n(\tb)$ does not vanish on the boundary of $\Delta$
but vanishes on $X_n$. It is easy to see that the cardinality of $X_n$ is 
$\dim \Pi_{n-1}^2$. Consequently,  the corresponding $\PP_n$ has $\dim
\Pi_{n-1}^2$ common zeros in $X_n^*$, so that the Gaussian cubature 
with respect to $w_{\frac{1}{2}}$ exists, which takes the form 
\eqref{GaussCuba}.  To obtain the formula for $\mu_j^{(n)}$, we note that 
the change of variables \eqref{x-y} shows that \eqref{GaussCuba} is the 
same as 
\begin{align*}
 \frac{1}{|\Delta |}  \int_{\Delta} f(\tb) \left[\TS_{1,1,-2}(\tb)\right]^2 d\tb 
   = \frac{2}{(n+2)^2} \sum_{j \in \Lambda_{n+2}^\circ} 
     \left[\TS_{1,1,-2}(\tfrac{\jb}{n+2})\right]^2
          f(\tfrac{\jb}{n+2}),
\end{align*}
where the last step follows from the fact that $\TS_{1,1,-2}(\tb)$ vanishes
on the boundary of $\Delta$. This is exactly the cubature \eqref{trigCuba}
applied to the function $f(\tb) \left[\TS_{1,1,-2}(\tb)\right]^2$. 
\end{proof}

We note that $w_{\frac{1}{2}}$ is only the second example of a weight function 
for which a Gaussian cubature exists for all $n$. Furthermore, the example
of $w_{-\frac{1}{2}}$ shows that the existence of Gaussian cubature does 
not hold for all weight functions in the family of $w_\alpha$. In \cite{SX}, a 
family of weight functions depending on three parameters that are on a region 
bounded by a parabola and two lines is shown to ensure Gaussian cubature 
for all permissible parameters. 

The set of nodes of a Gaussian cubature is poised for polynomial interpolation,
which shows that there is a unique polynomial $P$ in $\Pi_{n-1}^2$ such that
$P(\xi) = f(\xi)$ for all $\xi \in X_n^*$, where $f$ is an arbitrary function on 
$X_n^*$. In fact the general theory (cf.  \cite{DX}) immediately gives the 
following result: 

\begin{prop}
The unique interpolation polynomial of degree $n$ on $X_n^*$ is
given by
\begin{equation} \label{interp-tri-interior}
 \CL_n f(x) =  \sum_{j_1=1}^{n+1}\sum_{j_2=1}^{n+1- j_1} f(\tfrac{j}{n})
                       \mu_j^{(n)}K_n(\tb, \tfrac{j}{n}),
\end{equation}
where $x = (x_1,x_2)$ and $z(\tb) = x_1+ix_2$, and $\CL_n f$ is of total degree 
$n-1$ in $x$. 
\end{prop}

Under the mapping \eqref{x-y}, the above interpolation corresponds to a
trigonometric interpolation on $\Lambda_{n+2}^\circ$. Since 
$\mu_j^{(n)}= \frac{2}{(n+2)^2} |\TS_{1,1,-2}(\tfrac{\jb}{n})|^2$, it follows readily that 
$$
  \mu_j^{(n)}K_n(\tb, \tfrac{j}{n}) 
  = \frac{2}{(n+2)^2} \TS_{1,1,-2}(\tfrac{\jb}{n})
   \sum_{\kb \in \Lambda_{n+2}^\circ} \frac{\TS_{\kb}(\tb) \TS_{\kb}(\tfrac{\jb}{n})}
       {\TS_{1,1,-2}(\tb)}, 
$$ 
where $\jb = (j_1,j_2,-j_1-j_2)$. Consequently, the interpolation 
\eqref{interp-tri-interior} is closely related to the trigonometric 
interpolation on the triangle given in Theorem \ref{prop:1st-interpo-Dela}.  

\subsection{Cubature formula and Chebyshev polynomials of the first kind}

As discussed in the previous subsection, Gaussian cubature does not exist 
for $w_{-\frac12}$ on $\Delta^*$. It turns out, however, that there is a 
Gauss-Lobatto type cubature for this weight function whose nodes
are the points in the set 
\begin{align*}
Y_n := \Lambda_n 
   = \{(\tfrac{k_1}{n}, \tfrac{k_2}{n}) | k_1 \ge 0, k_2 \ge 0, 
        k_1+k_2 \le n\},
\end{align*}
which is located inside $\Delta$. Again we let $Y_n^*$ denote the image of 
$Y_n$ in $\Delta^*$ under the mapping \eqref{x-y}. 

The cardinality of $Y_n$ is $\dim \Pi_n^2$. According to a characterization 
of cubature formula in the language of ideals and varieties (\cite{X00}), there will 
be a cubature formula of degree $2n-1$ based on $Y_n^*$ if there are $n+2$ 
linearly independent polynomials of degree $n+1$ that vanish on $Y_n^*$.
In other words, $Y_n^*$ is the variety of an ideal generated by $n+2$ 
algebraic polynomials of degree $n+1$, which are necessarily  quasi-orthogonal 
in the sense that they are orthogonal to all polynomials of degree $n-2$.  
Such an ideal and a basis are described in the following proposition. 

\begin{prop} \label{prop1}
The set $Y_n^*$ is the variety of the polynomial ideal 
\begin{align} \label{ideal}
  \langle T_0^{n+1} - T_1^n, \quad T_k^{n+1} - T_{k-1}^{n-1} \quad (1 \le k \le n),
   \quad T_{n+1}^{n+1} - T_{n-1}^n \rangle. 
\end{align}
\end{prop}

\begin{proof}
The mapping \eqref{x-y} allows us to work with $T_k^n(\tb)$ on $Y_n$. Using 
the explicit formula for $\TC_{k,n-k,-n}$ in \eqref{TC_cos} , we have 
\begin{align*}
  T_k^n(\tb) =&  \frac{1}{3} \left[ e^{ \frac{i \pi}{3} (n -2k)(t_2-t_3)} \cos n \pi t_1
     \right . \\
      & \left .+ 
  e^{ \frac{i \pi}{3} (n -2k)(t_3-t_1)} \cos n \pi t_2 + e^{ \frac{i \pi}{3} (n -2k)(t_1-t_2)}      
          \cos n \pi t_3 \right].
\end{align*}
Since $n +1 -2 k= n-1-2(k-1)$, the above formula implies that 
\begin{align*}
T_k^{n+1}(\tb) - T_{k-1}^{n-1}(\tb) =   \frac{-1}{3} & \left[  
   e^{ \frac{i \pi}{3} (n -2k)(t_2-t_3)}      \sin \pi t_1\sin m \pi t_1     \right . \\
         &\,  + 
  e^{ \frac{i \pi}{3} (n -2k)(t_3-t_1)} \sin \pi t_2 \sin n \pi t_2 \\
      &\, \left .  +
      e^{ \frac{i \pi}{3} (n -2k)(t_1-t_2)}      
          \sin \pi t_3 \sin n \pi t_3 \right],
\end{align*}
from which it follows immediately that $T_k^{n+1} - T_{k-1}^{n+1}$ vanishes
on $Y_n$ for $1 \le k \le n$. Furthermore, a tedious computation shows that 
\begin{align*}
&  T_k^n(\tfrac{j_1}{n}, \tfrac{j_2}{n}, \tfrac{-j_1-j_2}{n} ) =
   \frac{-1}{6}  e^{ \frac{-2 i \pi}{3n} (n +1)j_1+j_2} \left [
      (-2 + e^{2 i \pi j_2} + e^{2 i \pi(j_1+j_2)} )      \right . \\
      & \left .+ 
     e^{2 i \pi \frac{j_2}{n}} (-2  e^{2 i \pi j_2} + e^{2 i \pi(j_1+j_2)}  +1)
    +      e^{2 i \pi \frac{j_1+j_2}{n}} (-2
               e^{2 i \pi j_1} +e^{2 i \pi j_2 } +1)\right],  
\end{align*}
which is evidently zero whenever $j_1$ and $j_2$ are integers. Finally, we note
that $T_{n+1}^{n+1} - T_{n-1}^n$ is the conjugate of $T_0^{n+1} - T_1^n$. 
This completes the proof.
\end{proof}

As mentioned before, this proposition already implies that a cubature formula
of degree $2n-1$ based on the nodes of $Y_n^*$ exists. 

\begin{thm}
For the weight function $w_{-\frac{1}{2}}$ on $\Delta^*$ the cubature 
formula 
\begin{equation} \label{GaussCuba2}
  c_{- \frac{1}{2}} \int_{\Delta^*} f(x,y) w_{-\frac{1}{2}}(x,y) dxdy = 
     \frac{1}{3n^2}  \sum_{k_1=0}^{n}\sum_{k_2=0}^{k_1}
         \lambda_{k}^{(n)} f(\tfrac{k_1}{n}, \tfrac{k_2}{n}), \quad \forall f \in \Pi_{2n-1}, 
\end{equation}
holds, where $\lambda_k^{(n)}= \lambda_{k,n-k,-n}^{(n)}$ with 
$\lambda_\jb^{(n)}$ given by \eqref{lambda_j}. 
\end{thm} 


In fact, the change of variables back to $\Delta$ shows that the cubature
\eqref{GaussCuba2} is exactly of the form \eqref{trigCuba}, which gives the
values of $\lambda_k^{(n)}$. In particular, this gives a proof of Theorem  
\ref{thm:trigCuba}. 

Another way of determining the cubature weight $\lambda_k^{(n)}$ is to work
with the reproducing kernel and the polynomial ideal. This approach also
gives another way of determining the formula for polynomial interpolation.
Furthermore, it fits into a general framework (see \cite{X97} and the references
therein), and the intermediate results are of independent interest.  In the rest 
of this section, we follow through with this approach. 


The recursive relation \eqref{recurT} for $T_k^m$ can be conveniently written 
in matrix form. Recall that for $T_k^m$, the orthonormal polynomials are
$$
 P_0^m = \sqrt{3} T_0^m, \quad P_k^m= \sqrt{6}T_k^m \quad (1 \le k \le m-1),
  \quad P_m^m=\sqrt{3}T_m^m
$$
for $m >0$ and $P_0^0=1$. In the following we treat $\PP_m$ defined in 
\eqref{PP_n} as a column vector. The relation \eqref{conjugT}
shows an important relation
\begin{equation}\label{JPP}
   \overline{ \PP_m } = J_{m+1}\PP_m, \qquad J_m = \left[ \begin{matrix} 
         \bigcirc  &  & 1 \\
            &  \bddots &  \\
        1 &  & \bigcirc \\
         \end{matrix} \right]. 
\end{equation}
Using these notations, we have the following three-term relations:

\begin{prop} \label{prop2}
Define $\PP_{-1} =0$.  For $m \ge 0$, 
\begin{align} \label{3term}  
\begin{split}
  z \PP_m & = A_m \PP_{m+1} + B_m \PP_m + C_m \PP_{m-1}, \\
 \bar z  \PP_m & = C_{m+1}^\tr \PP_{m+1} + B_m^\tr \PP_m + 
      A_{m-1}^\tr \PP_{m-1},
\end{split}
\end{align}
where $A_m$, $B_m$ and $C_m$ are given by 
$$
  A_m= \frac{1}{3} \,\left[ \begin{matrix} 
        1 &  &  & \bigcirc  & 0\\
       &  \ddots  & & & \vdots\\
       & &  1 & & \\
   \bigcirc & &  & \sqrt{2} & 0 
\end{matrix} \right]\;,\;\;
B_m=\frac{1}{3}\,\left[ \begin{matrix}
 0 & \sqrt{2} & 0 & \dots & &\bigcirc\\
  0 &  & 1 &  &     &\\
\vdots &  &  & \ddots& &  \\
   0&  \bigcirc& & & 1& 0 \\
    0  & 0 & \ldots & &0 &\sqrt{2}\\
\end{matrix} \right]\;.
$$
and $C_m = J_{m} A_{m-1}^\tr J_{m-1}$. 
\end{prop} 

In fact, the first relation is exactly \eqref{recurT}. To verify the second 
relation, take the complex conjugate of the first equation and then 
use \eqref{JPP}. 

In the following we slightly abuse the notation by writing $\PP_n(\tb)$ in 
place of $\PP_n(z)$, when $z=z(\tb)$ is given by \eqref{z}. Similarly, we write
$K_n(\tb,\sb)$ in place of $K_n(z, w)$ when $z = z(\tb)$ and $w = w(\sb)$. 
Just like the case of real orthogonal polynomials, the three-term relation 
implies a Christoph-Darboux type formula (\cite{DX})

\begin{prop}\label{prop3}
For $m \ge 0$, we have 
\begin{equation}\label{C-D}
  (z -w) K_n(z,w) := \PP_{n+1}^\tr(z) A_n^\tr\, \overline{\PP_n(w)} - 
           \PP_n^\tr(z) C_{n+1}^\tr \overline{\PP_{n+1}(w)}. 
\end{equation}
\end{prop}


Next we incorporate the information on the ideal in \eqref{ideal} into the kernel 
function. The $n+2$ polynomials in the ideal can be written as a set  
\begin{equation} \label{QQ}
     \QQ_{n+1}  = \PP_{n+1}  -  \Gamma_0\PP_n  -  \Gamma_1\PP_{n-1},
\end{equation}
which we also treat as a column vector below, where 
$$
    \Gamma_0 = \left[ \begin{matrix} 0& \frac{1}{\sqrt{2}} &0 &\ldots &0 \\
                    \bigcirc  &\cdots  & \bigcirc &\cdots  & \bigcirc \\
                   0& \ldots & 0 & \frac{1}{\sqrt{2}}& 0
        \end{matrix} \right] \; \quad\hbox{and}\quad  
            \Gamma_1 = \left[ \begin{matrix} 0& \ldots &0&\ldots & 0 \\
                    \sqrt{2} &   &  & & \bigcirc\\ 
                                & 1 &  & & \\
                   &  &\ddots  & &   \\
                    &   &  & 1 & \\
                     \bigcirc & & & &\sqrt{2}\\                 
                   0& \ldots & 0 & \ldots& 0
        \end{matrix} \right] \; \quad  
$$

\begin{prop}
For $n \ge 0$, define 
$$
 \Phi_n(\tb,\sb) := \frac{1}{2} \left[K_n(\tb,\sb) +K_{n-1}(\tb,\sb) \right] - \frac{1}{2}
   \left[  T_0^n(\tb) \overline{T_0^n(\sb)} + \overline{T_0^n(\tb)} T_0^n(\sb) \right],
$$
and we again write $\Phi_n(z,w)$ in place of  $\Phi_n(\tb,\sb)$ for 
$z = z(\tb)$ and $w=w(\sb)$. Then 
\begin{equation} \label{Phi-QQ}
  (z-w) \Phi_n(z,w) = \QQ_{n+1}^\tr(z) A_n^\tr M_n\overline{\PP_n(w)} - 
           \PP_n^\tr(z)  M_n C_{n+1}^\tr \overline{\QQ_{n+1}(w)} 
\end{equation}
where $M_n = \diag\{-1/3,1/2,\ldots,1/2,-1/3\}$ is a diagonal matrix. 
\end{prop}

\begin{proof}
Inserting $\PP_{n+1} = \QQ_{n+1} + \Gamma_0 \PP_n+ \Gamma_1\PP_{n-1}$
into \eqref{C-D} gives 
\begin{align*}
  (z-w) K_n(z,w) = \, & \QQ_{n+1}^\tr(z) A_n^\tr \,\overline{\PP_n(w)} - 
         \PP_n^\tr(z) C_{n+1}^\tr \overline{\QQ_{n+1}(w)} \\
           & + \PP_{n}^\tr(z) \Gamma_0^\tr A_n^\tr \,\overline{\PP_n(w)} - 
           \PP_n^\tr(z) C_{n+1}^\tr\Gamma_0 \overline{\PP_{n}(w)} \\
           & + \PP_{n-1}^\tr(z) \Gamma_1^\tr A_n^\tr \,\overline{\PP_n(w)} - 
           \PP_n^\tr(z) C_{n+1}^\tr\Gamma_1 \overline{\PP_{n-1}(w)}.
\end{align*}
Let $D_n$ be a $n\times n$  matrix. Then similarly we have
\begin{align*}
  & (z-w) \PP_n^\tr(z)D_n\overline{\PP_n(w)} = 
    \,  \QQ_{n+1}^\tr(z) A_n^\tr D_n \overline{\PP_n(w)} - 
         \PP_n^\tr(z) D_n C_{n+1}^\tr \overline{\QQ_{n+1}(w)} \\
       & \quad + \PP_{n}^\tr(z) (A_n\Gamma_0+B_n)^\tr D_n \overline{\PP_n(w)} - 
           \PP_n^\tr(z)D_n( C_{n+1}^\tr\Gamma_0 +B_n^\tr)\overline{\PP_{n}(w)} \\
       & \quad+ \PP_{n-1}^\tr(z) (A_n\Gamma_1 +C_n)^\tr D_n \overline{\PP_n(w)} - 
           \PP_n^\tr(z) D_n(C_{n+1}^\tr\Gamma_1 +A_{n-1}^\tr)\overline{\PP_{n-1}(w)}.
\end{align*}
In order to prove the stated result, we take the difference of the two equations
and choose $D_n$ such that 
\begin{align} \label{TwoCond} 
\begin{split}
 & (A_n\Gamma_0+B_n)^\tr D_n - D_n(C_{n+1}^\tr\Gamma_0 +B_n^\tr) 
      = \Gamma_0^\tr A_n^\tr - C_{n+1}^\tr\Gamma_1, \\
 &  (A_n\Gamma_1 +C_n)^\tr D_n = \Gamma_1^\tr A_n^\tr,\qquad 
    D_n( C_{n+1}^\tr\Gamma_1 +A_{n-1}^\tr)= C_{n+1}^\tr \Gamma_1. 
\end{split}
\end{align}
It turns out that the unique solution of the above three equations is given 
by $D_n = \diag\{2/3,1/2,1/2, \ldots, 1/2, 2/3\}$, which is completely determined
by the last two equations and then verified to satisfy the first equation. 
As a consequence of this choice of $D_n$, we obtain 
\begin{align*}
 (z-w)\left[K_n(z,w) - \PP_n^\tr(z)D_n\overline{\PP(w)}\right]
  & =  \QQ_{n+1}^\tr(z) A_n^\tr (I- D_n) \overline{\PP_n(w)} \\
      & - \PP_n^\tr(z) (I-D_n) C_{n+1}^\tr  \overline{\QQ_{n+1}(w)}
\end{align*} 
where $I$ denotes the identity matrix. Setting $M_n = I -D_n$ we see that 
the right hand side agrees with the right hand side of \eqref{Phi-QQ}. 
Finally we note that 
\begin{align*}
& \PP_n^\tr(z)D_n \overline{\PP_n(w)} =\frac{1}{2} \PP_n^\tr(z) \overline{\PP_n(w)}
 + \left[P_0^n(z)\overline{P_0^n(w)}  + P_n^n(z)\overline{P_n^n(w)} 
  \right]    \\
&\qquad \qquad = \frac{1}{2} \left[ K_n(z,w) - K_{n-1}(z,w)\right]
+ \frac{1}{2} \left[T_0^n(z)\overline{T_0^n(w)}  + \overline{T_0^n(z)}T_0^n(z)\right],    
\end{align*}
which completes the proof. 
\end{proof}

As a consequence of the above proposition, it is readily seen that 
$ \Phi_n( \tfrac{j}{n}, \tfrac{k}{n} ) = 0$ for $ j \ne k$. 
Furthermore, we note that $\Phi_n(\tb,\sb)$ is in fact a real valued function. 
Consequently, we have proved the following.

\begin{prop}
The unique interpolation polynomial on $Y_n^*$ is given by
$$ 
 \CL_n^* f(x,y) =  \sum_{j_1=0}^n\sum_{j_2=0}^{n-j_1} f(\tfrac{j}{n})
        \frac{\Phi_n(\tb, \tfrac{j}{n})}{\Phi_n(\tfrac{j}{n},\tfrac{j}{n})}
$$
where $z(\tb) = x+iy$, and $\CL_n^* f$ is of total degree $n$ in $x$ and $y$. 
\end{prop}

In fact, upon changing variables \eqref{z}, this interpolation polynomial is exactly 
the trigonometric interpolation function $\CL_n^*f (\tb)$ at \eqref{interpo-Dela}. 
Here the uniqueness of the interpolation follows from the general theory. 
Moreover, integrating $\CL_n^* f$ with respect to $w_{-\frac{1}{2}}$ 
over $\Delta^*$ yields the cubature formula \eqref{GaussCuba2} with
\begin{equation} \label{eq:lambda}
\lambda_j^{(n)} = \left[  \Phi_n(\tfrac{j}{n}, \tfrac{j}{n}) \right]^{-1},
\end{equation}
which can be used to determine $\lambda_j^{(n)}$ if needed. The function 
$K_n(\tb,\sb)$, thus $\Phi_n(\tb,\sb)$, enjoys a compact formula, which 
facilitates the computation of $\Phi_n(\tfrac{j}{n}, \tfrac{j}{n})$. The
compact formula follows from the following identity. 

\begin{prop}
Let $D_n$ be defined in \eqref{Gamma}. Then 
$$
 K_n(\sb,\tb)= \sum_{k=0}^n \sum_{j=0}^k \wt T_j^k(\tb) \overline{\wt T_j^k(\sb)} 
     = \CP^+_t D_n(\tb-\sb).
$$
\end{prop}

The identity can be established from the definition of $T_j^k$ and $D_n$ as
in the proof of Theorem \ref{prop:1st-interpo-Dela} upon using the fact that 
$\CP_\tb^+ \CP_\sb^+ F(\tb - \sb) = \CP^+ F(\tb- \sb)$. We omit the details. From
this identity and \eqref{eq:lambda}, we immediately have a compact formula 
for the interpolation polynomial, whose trigonometric counterpart under the 
change of variables \eqref{x-y} is exactly the function $\ell_{\jb,n}^\Delta$
in Theorem  \ref{prop:interpo-Dela}. 

\medskip\noindent
{\bf Acknowledgment.} The authors thank the referees for their meticulous 
and prudent comments.


\begin{thebibliography}{99}

\bibitem{AB}
       V. V. Arestov and E. E. Berdysheva, 
       Tur\'an's problem for positive definite functions with supports in a hexagon,
       {\it  Proc. Steklov Inst. Math. 2001, Approximation Theory. Asymptotical 
        Expansions}, suppl. 1, S20--S29.

\bibitem{CS}
        J. H. Conway and N. J. A. Sloane, 
        \textit{Sphere Packings, Lattices and Groups},  3rd ed. 
        Springer, New York, 1999.

\bibitem{DM} 
        D. E. Dudgeon and R. M. Mersereau, 
        \textit{Multidimensional Digital Signal Processing}, 
        Prentice-Hall Inc, Englewood Cliffs, New Jersey, 1984.

\bibitem{DX}
        C. F. Dunkl and Yuan Xu, 
        \textit{Orthogonal polynomials of several variables}, 
        Encyclopedia of Mathematics and its Applications, vol. {\bf 81},
        Cambridge Univ. Press, 2001. 

\bibitem{E}
        W. Ebeling, 
        \textit{Lattices and Codes},
       Vieweg, Braunschweig/Wiesbaden, 1994. 
         
\bibitem{F}
        B.  Fuglede, 
        Commuting self-adjoint partial differential operators 
        and a group theoretic problem,  
        \textit{J. Functional Anal.} \textbf{16} (1974), 101-121.
        
\bibitem{Hi} 
        J. R. Higgins, 
        \textit{Sampling theory in Fourier and Signal Analysis, Foundations}, 
       Oxford Science Publications, New York, 1996.

\bibitem{Ko}
         M. Kolountzakis, 
         The study of translation tiling with Fourier analysis, 
          \textit{Fourier analysis and convexity}, 131--187, 
          Appl. Numer. Harmon. Anal., BirkhŠuser Boston, Boston,
          MA, 2004.
                 
\bibitem{K}
        T. Koornwinder, 
        Orthogonal polynomials in two varaibles which are eigenfunctions 
        of two algebraically independent partial differential operators,
        \textit{Nederl. Acad. Wetensch. Proc. Ser. A77} = \textit{Indag. Math}.
        \textbf{36} (1974), 357-381.      

\bibitem{K2}
        T. Koornwinder, 
        Two-variable analogues of the classical orthogonal polynomials, in
         \textit{Theory and applications of special functions}, 435--495, ed.
        R.~A. Askey, Academic Press, New York, 1975.

\bibitem{Ma}
        R. J. Marks II, 
        \textit{Introduction to Shannon Sampling and Interpolation Theory},  
        Springer-Verlag, New York, 1991. 

\bibitem{Mc1}
        B. J. McCartin,  
        Eigenstructure of the equilateral triangle, Part I: 
        The Dirichlet problem,           
        \textit{SIAM Review}, \textbf{45} (2003), 267-287. 

\bibitem{Mc2}
        B. J. McCartin,  
        Eigenstructure of the equilateral triangle, Part II: 
        The Neumann problem,           
        \textit{Math. Prob. Engineering}, \textbf{8} (2002), 517-539. 

\bibitem{My} 
        I. P.  Mysoviskikh,
         \textit{Interpolatory cubature formulas},  Nauka, Moscow, 1981.
        
\bibitem{P1}
        M. A. Pinsky,
         The eigenvalues of an equilateral triangle, 
         \textit{SIAM J. Math. Anal.}, \textbf{11} (1980), 819-827.
         
\bibitem{P2}
         M. A. Pinsky,
         Completeness of the eigenfunctions of the equilateral triangle, 
         \textit{SIAM J. Math. Anal.}, \textbf{16} (1985), 848-851.
                  
\bibitem{R}
          T. J. Rivlin, 
           \textit{An introduction to the approximation of functions},
           Dover Publ., New York, 1981.
 
\bibitem{SX}
        H. J. Schmid and Yuan Xu, 
        On bivariable Gaussian cubature formulae,
        \textit{ Proc. Amer. Math. Soc.} \textbf{122} (1994), 833-842. 
               
\bibitem{SO}
        I. H. Sloan and T. R. Osborn,
        Multiple integration over bounded and unbounded regions,  
        \textit{J. Comp. Applied Math.} \textbf{17} (1987) 181-196. 
        
\bibitem{St}
        A. Stroud, 
        \textit{Approximate calculation of multiple integrals},
        Prentice-Hall, Englewood Cliffs, NJ, 1971.

\bibitem{Sun}
       J. Sun, 
       Multivariate Fourier series over a class of non tensor-product 
       partition domains, 
       \textit{J. Comput. Math.} \textbf{21} (2003), 53-62.       
       
\bibitem{LS}
        J. Sun and H. Li,
        Generalized Fourier transform on an arbitrary triangular domain,
        \textit{Adv. Comp. Math.}, \textbf{22} (2005), 223-248. 
       
\bibitem{Sz}
       G. Szeg\H{o},
       \textit{Orthogonal Polynomials},
       Amer. Math. Soc. Colloq. Publ. Vol.23, Providence, 4th edition,
       1975.

\bibitem{X97}
         Yuan Xu,
         On orthogonal polynomials in several variables, 
         \textit{Special Functions, $q$-series and Related Topics}, 
         The Fields Institute for Research in Mathematical Sciences, 
         Communications Series, Volume \textbf{14}, 1997, p. 247-270.
         
\bibitem{X00} 
          Yuan Xu,
          Polynomial interpolation in several variables, cubature formulae,
          and ideals, \textit{Advances in Comp. Math.}, 
          \textbf{12} (2000), 363--376.

\bibitem{Z}
        A. Zygmund,
        \textit{Trigonometric series},
        Cambridge Univ. Press, Cambridge, 1959. 
        
\end{thebibliography}
\end{document}